\documentclass[11pt]{smfart}

\usepackage{amsmath,a4wide}
\usepackage{amssymb}
\usepackage{xspace}
\usepackage{euscript}
\usepackage{mathrsfs} 
\usepackage[all]{xy}
\usepackage{hyperref}
\usepackage[francais]{babel}
\usepackage[T1]{fontenc}
\usepackage{textcomp}
\usepackage{smfthm}

\newtheorem{lemme2}[subsubsection]{Lemme}	

\newtheorem{proposition}[subsection]{Proposition}	
\newtheorem{proposition2}[subsubsection]{Proposition}

\newtheorem{theoreme}[subsection]{Th\'eor\`eme}

\newtheorem{corollaire}[subsection]{Corollaire}

\theoremstyle{definition}
\newtheorem{definition}[subsection]{D\'efinition}
\newtheorem{definition2}[subsubsection]{D\'efinition}

\newtheorem{remarque}[subsection]{Remarque}
\newtheorem{remarque2}[subsubsection]{Remarque}

\newenvironment{demo}{\noindent{\textit{D\'emonstration. --- }}}{~\qedsymbol \vspace{4mm}}

\renewcommand{\theenumi}{\textbf{\roman{enumi}}} 

\renewcommand\theequation{\thesection.\textbf{\alph{equation}}} 

\makeatletter
\@addtoreset{equation}{subsection}
\makeatother

\renewcommand\thesubsubsection{\textbf{\thesubsection.\arabic{subsubsection}}}
\renewcommand\thesubsection{\textbf{\thesection.\arabic{subsection}}}

\newcommand{\Z}{\mathbb Z}
\newcommand{\Q}{\mathbb Q}
\newcommand{\R}{\mathbb R}
\newcommand{\C}{\mathbb C}
\newcommand{\F}{\mathbb F}

\newcommand{\Gm}{{\mathbb G}_m}

\newcommand{\Hom}{\underline{\mathrm{Hom}}}

\newcommand{\Spec}{\mathrm{Spec}}

\newcommand{\iso}{\buildrel{\sim}\over\rightarrow}

\renewcommand{\L}{\mathrm{L}}		

\def\commutatif{\ar@{}[rd]|{\circlearrowleft}}
\def\cartesien{\ar@{}[rd]|{\square}}

\renewcommand{\lim}{{\mathrm{lim}}} 

\def\egazéro#1#2{[{\bf \'EGA}~$0_{\textsc{#1}}$~#2]}		

\numberwithin{equation}{subsubsection}
\renewcommand{\theequation}{\thesubsection.\Alph{equation}}

\newcommand{\isolong}{\buildrel{\sim}\over\longrightarrow}
\newcommand{\gsp}{\mathrm{GSp}}
\newcommand{\acal}{\mathcal{A}}
\newcommand{\bcal}{\mathcal{B}}
\newcommand{\mcal}{\mathcal{M}}
\newcommand{\mcalb}{\bar{\mathcal{M}}}
\newcommand{\acalb}{\bar{\mathcal{A}}}
\newcommand{\Qlb}{\bar{\Q}_\ell}
\newcommand{\pH}{{}^p\mathcal{H}}
\newcommand{\RPsi}{\mathrm{R}\Psi}
\newcommand{\IC}{\underline{\mathcal{IC}}}

\title{Sur une conjecture de Kottwitz au bord}
\author{Beno\^it Stroh}
\date{15 janvier 2011}
\email{stroh@math.univ-paris13.fr}
\address{C.N.R.S, Université Paris 13, LAGA,
99 avenue J.B. Clément,
93430 Villetaneuse
France}

\begin{document}

\begin{abstract}
Nous nous intéressons à la cohomologie d'intersection de la compactification minimale des variétés de Siegel à certaines places de mauvaise réduction. 
Nous calculons la trace semi-simple du morphisme de Frobenius sur les fibres des cycles proches du complexe d'intersection. Nous obtenons une généralisation commune de résultats de Morel et de Haines et Ngô.
\end{abstract}

\begin{altabstract}
We are interested in the intersection cohomology of the minimal compactification of Siegel modular varieties at some places of bad reduction.
We compute the semi-simple trace of the Frobenius morphism on the fibers of the nearby cycles of the intersection complex. We obtain a common generalization of results of Morel and Haines and Ngô.
\end{altabstract}

\maketitle

\tableofcontents

\section{Introduction}

Calculer les fonctions~$L$ des variétés de Shimura est une des clés de voûte du programme de Langlands.
Ces fonctions sont définies par des produits eulériens et se calculent donc place par place.
Soit~$p$ un nombre premier et~$X$ une variété de Shimura PEL de type A ou C. Le premier cas à considérer est celui où~$X$ a bonne réduction en~$p$.
Lorsque~$X$ est de plus supposée propre sur~$\Spec(\Z_p)$, le calcul du facteur~$L$ en~$p$ de la représentation galoisienne non ramifiée
$$\mathrm{H}^\bullet(X\times \Spec(\bar{\Q}_p), \Q_\ell)$$
de $\mathrm{Gal}(\bar{\Q}_p/\Q_p)$ (où~$\ell$ désigne n'importe quel nombre premier différent de~$p$) a été effectué par Kottwitz~\cite{Counting@Kottwitz} modulo certains lemmes fondamentaux prouvés depuis par Ngô et Waldspurger.

Lorsque $X$ est toujours supposé lisse sur~$\Spec(\Z_p)$ mais plus propre, la représentation galoisienne 
$$\mathrm{H}^\bullet(X\times \Spec(\bar{\Q}_p), \Q_\ell)$$
est non ramifiée mais se comporte de manière plutôt pathologique (elle n'est par exemple pas pure). Il est plus intéressant et plus naturel d'étudier plutôt la cohomologie d'intersection
$$\mathrm{IH}^\bullet(X^*\times \Spec(\bar{\Q}_p), \Q_\ell)$$
où~$X^\star$ est la compactification minimale de~$X$, construite par Faltings et Chai dans le cas des variétés de Siegel~\cite{Deg@FaltingsChai} et par Lan dans le cas PEL général~\cite{Compact@Lan}. Avant même de pouvoir mettre en oeuvre la stratégie de Kottwitz, on se heurte à un problème purement géométrique qui n'apparaissait pas dans le cas propre~: calculer la trace du morphisme de Frobenius sur les fibres du complexe d'intersection de~$X^\star$. Ce problème a été résolu par Morel dans~\cite{Unitaire@Morel} et~\cite{Poids@Morel} (les détails sont écrits dans le cas des variétés de Siegel et des variétés de Shimura associées aux groupes unitaires sur~$\Q$ mais la méthode est complètement générale). L'idée de Morel repose sur une nouvelle caractérisation du prolongement intermédiaire d'un faisceau pervers pur sur les variétés sur les corps finis~\cite[th.\ 3.1.4]{Poids@Morel}. 
Cette construction s'applique aux variétés de Shimura qui ont bonne réduction en~$p$ car le faisceau pervers décalé~$\Q_\ell$ est pur sur les variétés lisses sur~$\Spec(\F_p)$.

Morel a ensuite mené à bien la stratégie de Kottwitz~(\cite{Livre@Morel} et~\cite{Siegel@Morel}) et en a déduit une expression automorphe du facteur~$L$ en~$p$ des variétés de Siegel et des variétés de Shimura associées à un groupe unitaire sur~$\Q$.

Concentrons-nous à présent sur le cas de mauvaise réduction où~$X$ n'est plus supposée lisse sur~$\Spec(\Z_p)$.
Dans ce cas, la fonction~$L$ en~$p$ définie par Serre~\cite{Zeta@Serre} est trop compliquée à étudier car elle n'est plus directement reliée à la géométrie de la fibre spéciale sur~$\Spec(\F_p)$. Rapoport a défini~\cite{Bad@Rapoport} une modification de la fonction~$L$ en~$p$, la fonction~$L$ semi-simple en~$p$, qui se calcule via la formule de Lefschetz et est donc reliée à la géométrie de la fibre spéciale.

Supposons~$X$ propre et la mauvaise réduction de type parahorique. Contrairement au cas de bonne réduction propre, il y a un problème local à résoudre~: calculer la trace semi-simple du morphisme de Frobenius sur les fibres du complexe~$\RPsi(\Q_\ell)$ des cycles proches. La réponse, à savoir une fonction de Bernstein centrale dans une algèbre de Hecke parahorique, a été conjecturée par Kottwitz. Haines et Ngô ont démontré cette conjecture dans~\cite{Kottwitz@HainesNgo}.

Dans cet article, nous tentons de trouver une généralisation commune, appelée conjecture de Kottwitz au bord, des travaux de Morel et de Haines et Ngô. Nous nous plaçons donc dans le cadre suivant~: $X$ est une variété de Shimura PEL qui n'est pas propre et qui a mauvaise réduction de type parahorique en~$p$. Nous supposons de plus que~$X$ est une variété de Siegel (mais nos résultats valent dans une généralité plus grande, voir la remarque~\ref{remarque_generalisation_PEL}). Soit~$X^*$ la compactification minimale de~$X$ construite dans~\cite{Minimale@Stroh} sur~$\Spec(\Z_p)$. Nous sommes intéressés par le facteur~$L$ semi-simple en~$p$ de la cohomologie d'intersection
$$\mathrm{IH}^\bullet(X^*\times \Spec(\bar{\Q}_p), \Q_\ell)$$
et il nous faut calculer la trace semi-simple du morphisme de Frobenius  sur les fibres du complexe~$\RPsi_{X^*}\circ \IC_{X^*}(\Q_\ell)$ des cycles proches du complexe d'intersection de~$X^\star$. Cette question généralise simultanément celles étudiées par Morel (cas sans mauvaise réduction) et par Haines et Ngô (cas sans bord). Le corollaire~\ref{corollaire_trace_ss} donne une expression explicite de cette trace semi-simple en terme de théorie des groupes.
Cette expression fait intervenir le produit d'une fonction de Bernstein et d'une somme portant sur divers sous-groupe paraboliques du groupe définissant la variété de Shimura~$X$.

Nous commençons par montrer dans la partie~\ref{section_cycles_proches} que la cohomologie d'intersection de~$X^*\times \Spec(\Q_p)$ est égale à la cohomologie d'intersection de~$X^*\times \Spec(\F_p)$ à valeur dans les cycles proches~$\RPsi_X(\Q_\ell)$. Ceci résulte en fait de l'énoncé local de commutation $$\RPsi_{X^*}\circ \IC_{X^*}(\Q_\ell)\: = \: \IC_{X^*}\circ \RPsi_X(\Q_\ell)\: .$$
Il peut alors paraître envisageable d'appliquer le théorème général~\cite[3.1.4]{Poids@Morel} pour calculer la trace semi-simple du morphisme de Frobenius sur les fibres du complexe $\IC_{X^*}\circ \RPsi_X(\Q_\ell)$ sur~$X^*\times \Spec(\F_p)$. Mais les cycles proches~$\RPsi_X(\Q_\ell)$ ne sont  en général pas purs et la pureté est l'hypothèse primordiale de~\cite[th.\ 3.1.4]{Poids@Morel}. Nous décomposons donc ces cycles proches~$\RPsi_X(\Q_\ell)$ en faisceaux pervers simples dans le groupe de Grothendieck. On peut appliquer le théorème~\cite[3.1.4]{Poids@Morel} à chacun de ces constituents simples puisque les faisceaux pervers simples sont purs.
On obtient alors une expression pour la trace semi-simple du morphisme de Frobenius sur les fibres de~$\RPsi_X(\Q_\ell)$, mais cette expression est non explicite car elle fait intervenir les multiplicités inconnues des facteurs de Jordan-Hölder de~$\RPsi_X(\Q_\ell)$.
Ce problème est résolu dans la partie~\ref{section_conservation_multiplicite} où l'on montre que certaines de ces multiplicités sont égales à celles des facteurs de Jordan-Hölder de~$\RPsi_{X'}(\Q_\ell)$ avec~$X'$ une strate de bord de la compactification minimale~$X^*$ de~$X$.
Nous utilisons à plusieurs reprises l'existence de compactifications toroïdales de~$X$, qui ont été construites dans~\cite{Compact@Stroh} et résolvent certaines des singularités de la compactification minimale~$X^*$.

\medskip

Je remercie chaleureusement Sophie Morel pour de nombreux conseils prodigués lors d'un séjour à Harvard. Je remercie également Alain Genestier, Ulrich Görtz et Bao Chau Ngô de m'avoir aidé à clarifier certains points de cet article. Je remercie enfin le rapporteur pour une relecture attentive.

\section{Prolongement intermédiaire et troncation par le poids}

Cette partie est consacrée à de brefs rappels du théorème fondamental de~\cite{Poids@Morel}. La situation considérée ici n'a rien avec voir avec les variétés de Siegel, aussi adopterons nous des notations indépendantes du reste de l'article. Soient $a$ un entier relatif, $k$ un corps fini et~$X$ un schéma
séparé de type fini sur $\Spec(k)$.
Notons $D^b_m(X,\Qlb)$ la catégorie dérivée bornée des faisceaux $\ell$-adiques mixtes sur~$X$.
Morel a introduit dans~\cite[3.1]{Poids@Morel} un foncteur de \og troncation en poids inférieur à~$a$\fg
$$\mathrm{w}_{\leq a} \: : \:  D^b_m(X,\Qlb) \: \longrightarrow \: D^b_m(X,\Qlb)$$
et un foncteur de \og troncation en poids strictement supérieur à~$a$\fg 
$$\mathrm{w}_{> a} \: : \:  D^b_m(X,\Qlb) \: \longrightarrow \: D^b_m(X,\Qlb)\: .$$
Cette paire de foncteurs de troncation est associé à une t-structure sur $D^b_m(X,\Qlb)$ notée $({}^\mathrm{w} D^{\leq a},{}^\mathrm{w} D^{> a})$.
Dans toute la suite de l'article, le terme de faisceau pervers sera relatif à la t-structure perverse autoduale habituelle et le foncteur cohomologique associé sera noté $\pH^0$.
Lorsque~$K$ est pervers mixte sur $X$, le complexe $\mathrm{w}_{\leq a} K$ est le plus grand sous-faisceau pervers de~$K$ de poids~$\leq a$ et $\mathrm{w}_{>a} K$ le plus grand quotient de poids~$>a$.
En particulier, si $K$ est pervers pur de poids~$a$ sur~$X$, on~a $\mathrm{w}_{\leq a} K=K$ et $\mathrm{w}_{> a} K=0$.
Lorsque~$K$ est un complexe mixte, le complexe $\mathrm{w}_{\leq a} K$ vérifie que $\pH^n(\mathrm{w}_{\leq a} K)$ est mixte de poids~$\leq a$ pour tout~$n\in \Z$ et $\mathrm{w}_{>a}K$ vérifie que $\pH^n(\mathrm{w}_{>a} K)$ est mixte de poids~$>a$.

\begin{remarque}\label{remarque_deux_notions_poids}
Si $K$ n'est pas pervers, le complexe $\mathrm{w}_{\leq a} K$ n'est en général pas mixte de poids~$\leq a$ au sens de Deligne. En effet, la condition garantissant que $\mathrm{w}_{\leq a} K$ soit mixte de poids~$\leq a$ est que $\pH^n(\mathrm{w}_{\leq a} K)$ soit mixte de poids $\leq n+a$ pour tout $n\in \Z$, et cette condition n'est pas vérifiée si $\mathrm{w}_{\leq a} K$ a des groupes de cohomologie perverse non nuls en degré~$<0$.
\end{remarque}

Soit $j:U\hookrightarrow X$ un ouvert non vide.
Le foncteur $\mathrm{w}_{\leq a}$ est relié au prolongement intermédiaire~$j_{!*}$ par le théorème suivant prouvé dans~\cite[3.1.4]{Poids@Morel}.

\begin{theoreme}[Morel] \label{theoreme_Morel} Soit $A$ un faisceau pervers pur de poids $a$ sur $U$. Il existe un isomorphisme canonique $j_{!*} A  \iso \mathrm{w}_{\leq a}\: j_* A$ où $j_*$ désigne le foncteur dérivé d'image directe.
\end{theoreme}

\begin{remarque} Une des propriétés surprenantes du foncteur de troncation~$\mathrm{w}_{\leq a}$ est son caractère triangulé~\cite[3.1.3.(iii)]{Poids@Morel}. Il définit en particulier un morphisme entre groupes de Grothendieck. Ce morphisme joue un rôle clé dans~\cite{Poids@Morel} comme dans cet article.
\end{remarque}

Pour les applications aux variétés de Shimura, il est plus commode de travailler avec une version stratifiée de ce théorème. Donnons-nous donc une stratification $X=\coprod_{k=1}^n S_k$ par des sous-schémas localement fermés. On suppose que $S_n=U$ et que l'adhérence de $S_k$ est égale à $\coprod_{l\leq k} S_l$ pour tout $k\leq n$. En particulier, $U$ est l'unique strate ouverte de~$X$ et~$S_0$ l'unique strate fermée. Notons $i_k$ l'immersion de $S_k$ dans $X$ ; on a donc $i_n=j$. Notons $\mathcal{T}_k$ le foncteur $i_{k\: *}\circ \mathrm{w}_{>a} \circ i_k^*$ et, pour tout sous-ensemble ordonné $I=\{n_1 \leq \cdots \leq n_r\}$ de $\{0,\cdots,n-1\}$ de cardinal $\sharp I=r$, notons $\mathcal{T}_I$ le composé $\mathcal{T}_{n_1}\circ \cdots \circ \mathcal{T}_{n_r}$. La proposition suivante est démontrée dans~\cite[3.3.5]{Poids@Morel}.

\begin{proposition} \label{proposition_Morel_alternee} Soit $A$ un faisceau pervers pur de poids $a$ sur $U$. Dans le groupe de Grothendieck de $D^b_m(X,\Qlb)$ on a l'égalité
$$j_{!*} \: A \: = \: \sum_{I} \: (-1)^{\sharp I} \cdot \mathcal{T}_I\left( j_*\:  A \right)$$
où $I$ parcourt les sous-ensembles ordonnés de $\{0,\cdots,n-1\}$.
\end{proposition}

\section{Compactifications des variétés de Siegel}

\subsection{Variétés de Siegel}

Soit $g$ un entier positif, $p$ un nombre premier, $n\geq 3$ un entier non divisible par~$p$ et $\mathcal{D}=\{d_1<d_2<\cdots<d_s\}$ un sous-ensemble ordonné de $\{1,\cdots,g\}$.
Soit $V=\oplus_{j=1}^{2g}\: \Z \cdot x_j$ le~$\Z$-module libre de rang~$2g$ de base $(x_j)$. On le munit de la forme symplectique de matrice 
$$J=\left( \begin{array}{cc} 0 & J' \\
 -J' & 0           \end{array} \right)$$
où~$J'$ est la matrice anti-diagonale de coefficients non nuls égaux à~$1$.
Notons
$$\acal_V \: \longrightarrow \: \Spec(\Z[1/n])$$
l'espace de modules des schémas abéliens~$A$ principalement polarisés de genre~$g$ munis d'une similitude symplectique entre le groupes des points de $n$-torsion $A[n]$ et le $\Z/n\Z$-module symplectique $V/nV$, et d'une chaîne de groupes finis et plats
$ H_1 \subset H_2 \subset \cdots \subset H_s \subset A[p] $
telle que~$H_i$ soit de rang $p^{d_i}$ pour $0\leq i\leq s$ et que $H_s$ soit totalement isotrope pour l'accouplement de Weil.
On appelle cette structure de niveau en~$p$ une structure de niveau parahorique de type~$\mathcal{D}$.
L'espace de modules $\acal_V$ est représentable par un schéma quasi-projectif sur $\Spec(\Z[1/n])$. Ce schéma n'est ni propre ni lisse sur $\Spec(\Z[1/n])$, mais est lisse sur~$\Spec(\Z[1/np])$.
Sa dimension relative notée $d_V$ est égale à $g(g+1)/2$.

\subsection{Compactification minimale} \label{subsection_compactification_minimale} D\'ecrivons une compactification canonique de $\acal_V$ appel\'ee compactification minimale, de Satake ou bien encore de Baily-Borel, et construite dans~\cite{Minimale@Stroh}.
Il nous faut au préalable introduire quelques notations d'algèbre linéaire.

Pour tout sous-espace totalement isotrope $V'$ facteur direct  de $V$, notons $V'^\perp$ son orthogonal et $C(V/V'^\perp)$ le cône des formes quadratiques semi-définies positives à radical rationnel sur $(V/V'^\perp)\otimes \R$.
Notons $\mathfrak{C}_V$ l'ensemble des sous-espaces totalement isotropes facteurs directs de~$V$ et $\mathcal{C}_V$  le quotient  de l'union disjointe
$$\coprod_{V'\in \mathfrak{C}_V} C(V/V'^\perp)$$
par la relation d'équivalence engendrée par les inclusions $C(V/V''^\perp)\subset C(V/V'^\perp)$ pour tous sous-espaces totalement isotropes facteurs directs $V''\subset V'$.
L'ensemble $\mathfrak{C}_V$ est un {complexe conique}~\cite[p. 129]{Deg@FaltingsChai} muni d'une action du groupe~$\gsp(V)$ des similitudes symplectiques de~$V$. Soit $\mathbb{V}^\bullet$ le la chaîne de réseaux de $V$ tel que
$$\mathbb{V}^i  \:= \: \bigoplus_{j=1}^{d_i}\: \Z \cdot x_j \oplus \bigoplus_{j=d_i+1}^{2g}p\: \Z \cdot x_j$$
$$\mathrm{et} \:\:\: \: \: \mathbb{V}^{s+i}\: = \: \bigoplus_{j=1}^{2g-d_{s+1-i}} \:  \Z\cdot x_j \oplus \bigoplus_{j=2g - d_{s+1-i}+1}^{2g} p\:\Z$$
pour $1\leq i\leq s$.
Posons
$\Gamma_V = \mathrm{Stab}_{\gsp(V)} \left(\mathbb{V}^\bullet\right) \cap \mathrm{Ker}\left(\gsp(V)\rightarrow \gsp(V/nV)\right)\: .$
C'est le sous-groupe de congruence de $\gsp(V)$ naturellement associé au problème de modules définissant~$\acal_V$.

D'après~\cite[th.\ 3.9]{Minimale@Stroh}, il existe un schéma $\acal_V^*$ projectif de type fini sur~$\Spec(\Z[1/n])$ qui contient $\acal_V$ comme ouvert dense. Notons
$$j \: : \: \acal_V \: \hookrightarrow \: \acal_V^*$$
cette immersion ouverte. Le schéma $\acal_V^*$ est muni d'une stratification paramétrée par l'ensemble fini~$\mathfrak{C}_V/\Gamma_V$. La strate associée à $V' \in \mathfrak{C}_V$ est une variété de Siegel $\acal_{V'}$ associée au module symplectique~$V'^\perp/V'$, de type parahorique donné par l'image de $\mathbb{V}^\bullet$ dans $V'/V^{'\perp}$ et de niveau principal en~$n$. Les relations d'incidences entre strates sont duales des relations d'inclusions dans~$\mathfrak{C}_V$. L'unique strate ouverte correspond au sous-espace $\{0\} \in \mathfrak{C}_V$ et est égale à~$\acal_V$. Notons
$$i_{V'} \: : \:  \acal_{V'} \: \hookrightarrow \: \acal_V^*$$
l'inclusion localement fermée de la strate $\acal_{V'}$ de $\acal_V^*$. On a donc $i_{\{0\}}=j$.

\begin{remarque2} La situation est particulièrement simple lorsque $\mathcal{D}=\{1,\cdots,g\}$, qui est le cas d'une structure de niveau iwahorique. En effet, les structures de niveau des variétés de Siegel $\acal_{V'}$ sont alors toutes iwahoriques. Dans le cas général o\`u $\mathcal{D}\neq\{1,\cdots,g\}$, il se peut que certains des strates~$\acal_{V'}$ soient sans structure de niveau en~$p$. Elles sont alors lisses sur $\Spec(\Z[1/n])$.
\end{remarque2}

\begin{remarque2} Bien que $\acal_V$ soit lisse sur $\Spec(\Z[1/np])$, ce n'est jamais le cas de $\acal_V^*$ lorsque $g\geq 2$. En effet, la variété complexe sous-jacente à $\acal_V^*$ est isomorphe à la compactification de l'espace localement symétrique $\acal_V(\C)$ construite par Satake, et cette dernière n'est pas une variété complexe lisse. Le schéma $\acal_V$ cumule donc deux tares : il a des singularités \og verticales \fg qui proviennent du bord et sont déjà présentes en caractéristique nulle, et des singularités \og horizontales \fg qui proviennent de la mauvaise réduction de $\acal_V$ sur~$\Spec(\F_p)$.
Notons que nous ne chercherons pas à donner un sens mathématique à ce concept de verticalité et d'horizontalité. L'origine de ces termes sera claire pour le lecteur qui dessinera~$\Spec(\Z[1/n])$ comme une courbe horizontale sur laquelle se fibre $\acal_V^*$.
Notons aussi que les strates $\acal_{V'}$ ont en général elles-mêmes mauvaise réduction sur~$\Spec(\F_p)$. Une grande partie de l'article sera consacrée à la compréhension des interactions entre singularités horizontales et verticales.
\end{remarque2}

\subsection{Compactifications toroïdales} \label{subsection_compact_toro} Il existe une famille de \og résolutions des singularités verticales \fg de $\acal_V^*$.
Les éléments de cette famille s'appellent compactifications toroïdales et sont indexés par le choix d'une décomposition polyédrale $\gsp(V)$-admissible $\mathfrak{S}_V$ du complexe conique~$\mathcal{C}_V$~\cite[déf.\ 3.2.3.1]{Compact@Stroh}.
La compactification toroïdale associ\'ee à~$\mathfrak{S}$ est notée $\acalb_V^\mathfrak{S}$ et est construite dans~\cite{Compact@Stroh}. C'est \textit{a priori} seulement un espace algébrique, mais si $\mathfrak{S}$ vérifie une hypothèse de polarisation~\cite[déf.\ IV.2.4]{Deg@FaltingsChai}, ce que l'on peut toujours supposer quitte à raffiner, $\acalb_V^\mathfrak{S}$ est un schéma projectif de type fini sur $\Spec(\Z[1/n])$. De plus, $\acalb_V^\mathfrak{S}$ est lisse sur~$\Spec(\Z[1/np]$ pour peu que~$\mathfrak{S}$ soit choisi assez fin, ce que l'on supposera. Choisissons dans tout le reste de l'article une décomposition assez fine~$\mathfrak{S}$ et pour alléger les notations, notons~$\acalb_V$ la compactification toroïdale associée. 

Notons $J:\acal_V \hookrightarrow \acalb_V$ l'immersion ouverte canonique. D'après~\cite[th.\ 3.9]{Minimale@Stroh}, il existe un morphisme $\pi$ de la compactification toroïdale dans la minimale rendant commutatif le diagramme
$$\xymatrix{ & \acalb_V \ar[d]^{\pi} \\
\acal_V \ar[ru]^{J} \ar[r]^j & \acal_V^*
}$$
La stratification de~$\acal_V^*$ induit une stratification sur~$\acalb_V$ paramétrée par $\mathfrak{C}_V/\Gamma_V$. L'avantage des compactifications toroïdales est que le voisinage étale de chacune de ces strates peut aisément être décrit de la manière suivante.

Soit $V'\in \mathfrak{C}_V$. Dans~\cite[par.\ 1.4]{Compact@Stroh} sont décrits un schéma abélien $\bcal_{V'}\rightarrow \acal_{V'}$ isogène \`a une puissance du sch\'ema abélien universel sur~$\acal_{V'}$, un torseur $\mcal_{V'}\rightarrow \bcal_{V'}$ sous un tore et un plongement en fibrés toriques localement de type fini
$$\xymatrix{\mcal_{V'} \ar[r] \ar[rd] & \mcalb_{V'} \ar[d] \\
 & \bcal_{V'}
}$$
de fibre type le pongement torique $E_{V'}\subset \bar{E}_{V'}$ associé à la restriction de l'éventail~$\mathfrak{S}$ à $C(V/V^{'\perp})$. Le fibré en plongement torique $\mcalb_{V'}$ est muni d'une strate fermée canonique qui s'obtient \textit{via} la théorie des plongements toriques en restreignant l'éventail~$\mathfrak{S}$ à l'intérieur du cône~$C(V/V^{'\perp})$.
Notons $\Gamma_{{V'}}$ le stabilisateur de $V'$ dans~$\Gamma_V$ et~$\Gamma_{V'}^{l}$ l'intersection
$\Gamma_{V'}^{l} \: = \:  \Gamma_{{V'}}\cap \mathrm{GL}(V')$.
Le sous-groupe discret~$\Gamma_{V'}^{l}$  d'un groupe linéaire agit sans points fixes sur $\mcalb_{V'}$ au dessus de~$\bcal_{V'}$ et le quotient est de type fini sur~$\bcal_{V'}$. Par construction des compactifications toroïdales~\cite[par.\ 3.1.8]{Compact@Stroh}, il existe un isomorphisme entre l'hensélisé de $\acalb_V$ le long de la strate paramétrée par~$V'$ et le quotient par $\Gamma_{V'}^{l}$ de l'hensélisé de $\mcalb_{V'}$ le long de sa strate fermée.

\begin{remarque2} Ainsi, les singularités de $\acalb_V$ sont exactement les singularités des variétés de Siegel $\acal_{V'}$ pour $V'\in \mathfrak{C}_V$. Cela explique pourquoi $\acalb_V$ est lisse sur~$\Spec(\Z[1/np])$. Cela explique également pourquoi $\acalb_V$ résoud en un certain sens les singularités verticales de~$\acal_V^*$. Encore une fois, nous ne cherchons pas à donner de sens mathématique à cette notion de résolution de singularités verticales. Disons simplement que $\acalb_V \rightarrow \acal_V^*$ est une résolution des singularités sur $\Spec(\Z[1/np])$ au sens habituel.
\end{remarque2}

\begin{remarque2} Contrairement au cas des compactifications toroïdales, il n'existe aucune description géométrique explicite du voisinage de la strate~$\acal_{V'}$ de la compactification minimale~$\acal_V^*$ pour $V'\in \mathfrak{C}_V$.
\end{remarque2}

\subsection{Systèmes de coefficients} \label{subsection_systemes_coefficients}

Soit $\ell$ un nombre premier différent de~$p$. Il existe une mani\`ere canonique d'associer \`a toute repr\'esentation $R$ du groupe alg\'ebrique $\gsp(V\otimes\Q)$ (que l'on confondra désormais avec l'ensemble de ses points sur~$\Q$) un système local $\ell$-adique $\mathcal{F}_V(R)$ sur~$\acal_V\times \Spec(\Z[1/nl])$~\cite[p.\ 238]{Deg@FaltingsChai}. Ce système local est découpé par des correspondances alg\'ebriques dans la cohomologie de puissances du schéma abélien universel sur~$\acal_V\times \Spec(\Z[1/nl])$.
Par exemple, le système local $\mathcal{F}_V(\mathrm{Std}_V)$ associé à la représentation standard $\mathrm{Std}_V$ de $\gsp(V\otimes \Q)$ sur $V\otimes \Q$ est le module de Tate relatif du schéma abélien universel.

\begin{remarque2} Nous suivons ici les conventions de~\cite{Poids@Morel}, qui sont duales de~\cite{Deg@FaltingsChai}.
\end{remarque2}

\noindent Le foncteur $R\mapsto \mathcal{F}_V(R)$ commute à la dualité des faisceaux localement constants et des représentations, c'est à dire vérifie
$$\Hom_{\acal_V}(\mathcal{F}_V(R),\Q_l)\: = \: \mathcal{F}_V(\mathrm{Hom}(R,\Q))\: .$$
Il vérifie également $\mathcal{F}_V(R)(1)=\mathcal{F}_V(R\otimes c)$ où $(1)$ désigne le twist de Tate et $c$ est la représentation de dimension un de $\gsp(V\otimes \Q)$ associée au morphisme facteur de similitude de~$\gsp(V\otimes \Q)$ dans~$\Gm$.

Le foncteur exact $R\mapsto \mathcal{F}_V(R)$ se dérive en un foncteur~$\mathcal{F}_V$ de la catégorie dérivée bornée~$D^b(\gsp(V\otimes \Q))$ des représentations du groupe algébrique~$\gsp(V\otimes \Q)$ dans la catégorie dérivée~$D^b_c(\acal_V)$ des complexes $\ell$-adiques bornés constructibles sur $\acal_V\times \Spec(\Z[1/nl])$. Notons que comme le groupe~$\gsp(V\otimes \Q)$ est réductif, les objets de~$D^b(\gsp(V\otimes \Q))$ sont (non canoniquement) somme directe des décalés de leurs objets de cohomologie. Il en est donc de même des complexes~$\mathcal{F}_V(R)$ pour $R\in D^b(\gsp(V\otimes \Q))$. Disons que~$\Gm$ agit avec des poids~$\leq a$ sur un espace vectoriel~$R$ si tous les caractères intervenant dans l'action de~$\Gm$ sur~$R$ sont de la forme~$z\mapsto z^b$ avec $b\leq a$.

\begin{definition2}\label{definition_representation_pure}
Un complexe $R \in D^b(\gsp(V\otimes \Q))$ est mixte de poids $\leq a$ avec $a\in \Z$ si le centre~$\Gm$ de~$\gsp(V\otimes \Q)$ agit avec des poids~$\leq a+n$ sur le groupe de cohomologie~$\mathrm{H}^n(R)$ pour tout~$n\in \Z$. De même pour les notions de mixité en poids $\geq a$ et de pureté en poids~$a$.
\end{definition2}

La restriction à~$\acal_V\times \Spec(\F_p)$ du complexe~$\mathcal{F}_V(R)$ est toujours mixte, et l'on dispose du lemme suivant~\cite[5.6.6]{Baily@Pink}.

\begin{lemme2} \label{poids_rep_faisceau} Soit $R \in D^b(\gsp(V\otimes \Q))$ un complexe mixte de poids $\geq a$ (resp. mixte de poids~$\geq a$ ou pur de poids~$a$). La restriction à~$\acal_V\times \Spec(\F_p)$ du complexe~$\mathcal{F}_V(R)$ est mixte de poids ponctuels $\leq -a$ (resp. $\geq -a$ ou égaux à~$-a$).
\end{lemme2}

\begin{remarque2} \label{remarque_pas_pur} Comme $\acal_V\times \Spec(\F_p)$ n'est pas lisse sur~$\Spec(\F_p)$, le complexe~$\mathcal{F}_V(R)$ n'est pas nécessairement pur de poids~$a$ au sens de Deligne  lorsque~$R$ est pur de poids~$a$. Ceci devient par contre vrai après restriction à tout sous-schéma de~$\acal_V\times \Spec(\F_p)$ lisse sur~$\Spec(\F_p)$.
\end{remarque2}

\begin{definition} \label{definition_troncation_rep} Soit $a\in \Z$. Le foncteur $\mathrm{w}_{\leq a}$ de $D^b(\gsp(V\otimes \Q)$ dans~$D^b(\gsp(V\otimes \Q)$ \linebreak associe à tout complexe~$R$ l'unique facteur direct~$R'$ de~$R$ tel que le centre~$\Gm$ de $\gsp(V\otimes \Q)$ agisse avec des poids~$\leq a$ sur le groupe de cohomologie~$\mathrm{H}^n(R)$ pour tout~$n\in \Z$. De même pour~$\mathrm{w}_{>a}$. Les foncteurs~$\mathrm{w}_{\leq a}$ et~$\mathrm{w}_{> a}$ sont triangulés et passent donc au groupe de Grothendieck.
\end{definition}

\begin{remarque2} Le lecteur remarquera la différence entre le poids~$a+n$ intervenant dans la définition~\ref{definition_representation_pure} et le poids~$a$ intervenant dans la définition~\ref{definition_troncation_rep}. Ainsi, $\mathrm{w}_{\leq a}$ n'est pas mixte de poids~$\leq a$ si~$R$ a des groupes de cohomologie non nuls en degré~$<0$. Voir la remarque~\ref{remarque_deux_notions_poids} pour l'analogue dans le monde des faisceaux $\ell$-adiques.
\end{remarque2}

\section{Cycles proches et prolongement intermédiaire} \label{section_cycles_proches} Nous allons prouver plusieurs énoncés de commutation entre le foncteur des cycles proches et les foncteurs de prolongement par zéro et d'image directe, ceci sur les compactifications toroïdales comme sur la compactification minimale. Notons~$\bar{\eta}$ le spectre d'une clôture algébrique de~$\Q_p$. Rappelons que nous avons noté $J:\acal_V \hookrightarrow \acalb_V$ l'immersion ouverte de la variété de Siegel dans une de ses compactifications toroïdales. Pour alléger les notations, nous désignerons par $J_*$ le foncteur image directe dérivée de~$D^b_c(\acal_V)$ dans~$D^b_c(\acalb_V)$. La proposition suivante est réminiscente de~\cite[par.\ 7.1]{TaylorWiles@GenestierTilouine}

\begin{proposition} \label{proposition_commutation_toroidale} Pour tout $R\in D^b(\gsp(V\otimes Q))$, les morphismes d'adjonction
\begin{eqnarray*}
\RPsi \circ J_* \: \mathcal{F}_V(R) & \longrightarrow &  J_* \circ \RPsi \: \mathcal{F}_V(R) \\
J_! \circ \RPsi \: \mathcal{F}_V(R) & \longrightarrow & \RPsi\circ  J_* \: \mathcal{F}_V(R)
\end{eqnarray*}
sont des isomorphismes dans $D^b_c((\acalb_V \times \Spec(\bar{\F}_p))\times \bar{\eta})$.
\end{proposition}

\begin{demo} Par dualité de Poincaré, il suffit de traiter le cas de $J_*$.
Commençons par traiter le cas de la représentation constante~$R$, cas où~$\mathcal{F}_V(R)$ est le faisceau constant~$\Q_\ell$.
Voir si le morphisme d'adjonction est un isomorphisme se teste localement pour la topologie étale. Il suffit donc de prouver un énoncé analogue sur les descriptions explicites du voisinage de la strate de bord de~$\acalb_V$ paramétrée par~$V'\in \mathfrak{C}_V$. Soit $J':\mcal_{V'} \hookrightarrow \mcalb_{V'}$ l'immersion ouverte au-dessus de~$\bcal_{V'}$ décrite dans la partie~\ref{subsection_compact_toro}. Il suffit de prouver que le morphisme d'adjonction
$$\RPsi_{\mcalb_{V'}} \circ J'_*\:  (\Q_\ell) \: \longrightarrow \:  J'_* \circ \RPsi_{\mcal_{V'}} (\Q_\ell)$$
est un isomorphisme.
Quitte à localiser pour la topologie étale sur~$\bcal_{V'}$, on peut supposer que le fibré en plongement torique~$J':\mcal_{V'}\hookrightarrow \mcalb_{V'}$ sur~$\bcal_{V'}$ est trivial de fibre le plongement torique $J'':E_{V'} \hookrightarrow \bar{E}_{V'}$.
Il suffit donc de prouver que le morphisme d'adjonction
$$\RPsi_{\bcal_{V'}\times \bar{E}_{V'}} \circ (\mathrm{Id}_{\bcal_{V'}}\times J'')_*(\Q_\ell) \: \longrightarrow \:  (\mathrm{Id}_{\bcal_{V'}}\times J'')_* \circ \RPsi_{\bcal_{V'}\times E_{V'}} (\Q_\ell)$$
est un isomorphisme. D'après la formule de Künneth, on a un isomorphisme de foncteurs
$$(\mathrm{Id}_{\bcal_{V'}}\times J'')_* \: \isolong \: \mathrm{Id}_{\bcal_{V'}} \boxtimes J''_*$$
et d'après la formule de produit pour les cycles proches~\cite{SemiStable@Illusie}, on a de même
$$\RPsi_{\bcal_{V'}\times \bar{E}_{V'}}  \: \isolong \: \RPsi_{\bcal_{V'}} \boxtimes \RPsi_{\bar{E}_{V'}}\: .$$
Il suffit donc de prouver que
$$\RPsi_{\bar{E}_{V'}}\circ J''_* (\Q_\ell) \: \isolong \: J''_* \circ \RPsi_{E_{V'}} (\Q_\ell)$$
ce qui résulte de SGA~VII.2 2.1.9
car le complémentaire de~$E$ dans~$\bar{E}$ est un diviseur à croisements normaux d'après nos hypothèses combinatoires.

Dans le cas général où $R\in D^b(\gsp(V\otimes \Q))$, la preuve est tout à fait similaire une fois incorporée dans le formalisme l'existence de compactifications des variétés de Kuga-Sato. En effet, grâce à des arguments standards de pléthysmes~\cite[p.\ 235]{Deg@FaltingsChai} et le principe de Liebermann~\cite[p.\ 219]{Deg@FaltingsChai}, on se réduit à prouver que
$$\RPsi \circ J_* \circ g_* (\Q_\ell) \:  \isolong \:  J_* \circ \RPsi \circ g_* (\Q_\ell)$$
où $g:G^s\rightarrow \acal_V$ désigne la variété abélienne universelle~$G$ sur~$\acal_V$ à une certaine puissance \linebreak $s\in \mathbb{N}$. D'après~\cite[th.\ VI.1.1]{Deg@FaltingsChai}, $g:G^s\rightarrow \acal_V$ admet des compactifications explicites \linebreak $\bar{g}:\overline{G^s}\rightarrow \acalb_V$ qui proviennent par image inverse de compactifications sur les variétés de Siegel sans niveau en~$p$. Choisissons une de ces compactifications et notons $J''':G^s\hookrightarrow \overline{G^s}$ l'immersion ouverte au-dessus de $J:\acal_V\hookrightarrow \acalb_V$. D'après le théorème de changement de base propre, il suffit de prouver que
$$\RPsi_{\overline{G^s}} \circ J'''_* (\Q_\ell) \:  \isolong \:  J'''_* \circ \RPsi_{G^s} (\Q_\ell)\: .$$
On raisonne alors exactement comme dans le cas précédent où $\mathcal{F}_V(R)=\Q_\ell$ en utilisant la description locale sur~$\bcal_{V'}$ de~$G^s\subset \overline{G^s}$ comme produit de~$\bcal_{V'}$ par un plongement torique $F_{V'}\subset \bar{F}_{V'}$~\cite[th.\ VI.1.13.(iii)]{Deg@FaltingsChai}.
\end{demo}

\begin{remarque} Le lecteur se convaincra aisément qu'un énoncé de commutation des cycles proches avec les images directes par une immersion ouverte ne saurait être vrai en toute généralité. Cet énoncé est vrai dans le cas de~$\acal_V$ car la compactification~$\acalb_V$ admet localement une structure produit~$\bcal_{V'}\times \bar{E}_{V'}$, que la mauvaise réduction de~$\acalb_V$ est concentrée dans le facteur~$\bcal_{V'}$ et que l'immersion ouverte de~$\acal_V$ dans~$\acalb_V$ ne concerne que~$\bar{E}_{V'}$.
\end{remarque}

\noindent La proposition~\ref{proposition_commutation_toroidale} jointe à la propreté de~$\acalb_V$ implique notamment que
$$\mathrm{H}^i\left( \acal_V \times \Spec(\bar{\Q}_p) ,\: \Q_\ell \right) \: = \: \mathrm{H}^i\left( \acal_V \times \Spec(\bar{\F}_p) ,\: \RPsi(\Q_\ell) \right)$$
comme représentation de $\mathrm{Gal}(\bar{\Q}_p/\Q_p)$ pour tout~$i\in \mathbb{N}$, ce qui démontre une conjecture de Haines~\cite[10.3]{Clay@Haines}. Rappelons que l'on a noté~$j$ l'immersion ouverte de~$\acal_V$ dans la compactification minimale~$\acal_V^*$.

\begin{corollaire} \label{corollaire_commutation_minimale} Pour tout $R\in D^b(\gsp(V\otimes Q))$, les morphismes d'adjonction
\begin{eqnarray*}
\RPsi \circ j_* \: \mathcal{F}_V(R) & \longrightarrow &  j_* \circ \RPsi \: \mathcal{F}_V(R) \\
j_! \circ \RPsi \: \mathcal{F}_V(R) & \longrightarrow & \RPsi\circ  j_* \: \mathcal{F}_V(R)
\end{eqnarray*}
sont des isomorphismes dans $D^b_c((\acal_V^* \times \Spec(\bar{\F}_p))\times \bar{\eta})$.
\end{corollaire}

\begin{demo} Il suffit de combiner la proposition~\ref{proposition_commutation_toroidale} et le théorème de changement de base propre pour le morphisme~$\pi$ de la compactification toroïdale~$\acalb_V$ dans la minimale~$\acal_V^*$.
\end{demo}

\noindent Rappelons que d'après~\cite{SemiStable@Illusie}, le foncteur $\RPsi$ envoie faisceaux pervers sur faisceaux pervers.

\begin{corollaire} \label{corollaire_commutation_minimale} Notons $d_V=g(g+1)/2$ qui est la dimension relative de $\acal$ sur~$\Spec(\Z[1/n])$. Pour toute représentation~$R$ de~$\gsp(V\otimes Q)$, il existe un isomorphisme canonique
$$\RPsi \circ j_{!*} (\mathcal{F}_V(R)[d_V]) \: \isolong \: j_{!*} \circ \RPsi (\mathcal{F}_V(R)[d_V])$$
dans la catégorie des faisceaux pervers sur~$\acal_V^*\times \Spec(\bar{\F}_p)$ munis d'une action compatible de~$\mathrm{Gal}(\bar{\Q}_p/\Q_p)$.
\end{corollaire}

\begin{demo} De l'égalité $\RPsi \circ j_{*} (\mathcal{F}_V(R)[d_V]) =  j_{*} \circ \RPsi (\mathcal{F}_V(R)[d_V])$ on tire l'égalité ${}^p\mathcal{H}^0(\RPsi \circ j_{*} (\mathcal{F}_V(R)[d_V]))={}^p\mathcal{H}^0(j_{*} \circ \RPsi (\mathcal{F}_V(R)[d_V])$. Comme $\RPsi$ est $t$-exact pour la perversité autoduale, on~a ${}^p\mathcal{H}^0 \circ \RPsi = \RPsi \circ {}^p\mathcal{H}^0$. On en déduit $\RPsi \circ {}^p j^0_{*} (\mathcal{F}_V(R)[d_V]) =  {}^p j^0_{*} \circ \RPsi (\mathcal{F}_V(R)[d_V])$. De même, $\RPsi \circ {}^p j^0_{!} (\mathcal{F}_V(R)[d_V]) =  {}^p j^0_{!} \circ \RPsi (\mathcal{F}_V(R)[d_V])$. On aboutit au résultat en utilisant la définition du foncteur de prolongement intermédiaire~$j_{!*}$ comme l'image perverse de~${}^p j^0_{!}$ dans~${}^p j^0_{*}$.
\end{demo}

\begin{remarque} \label{remarque_referee} Nous avons  montré l'égalité des foncteurs qui à tout objet~$R$ de la catégorie semi-simple des représentations de~$\gsp(V\otimes Q)$ associent les faisceaux pervers avec action galoisienne $\RPsi \circ j_{!*} (\mathcal{F}_V(R)[d_V])$ et $j_{!*} \circ \RPsi (\mathcal{F}_V(R)[d_V])$. Ces deux foncteurs sont exacts et se dérivent trivialement. On obtient donc l'égalité entre deux foncteurs de $D^b(\gsp(V\otimes Q))$ dans $D^b_c((\acal_V^* \times \Spec(\bar{\F}_p))\times \bar{\eta})$ d'après~\cite{Pervers@Beilinson}. Nous réutiliserons ce principe à plusieurs reprises : dès qu'une égalité est démontrée entre faisceaux pervers dépendant fonctoriellement d'une représentation de~$\gsp(V\otimes Q)$, l'égalité analogue est  valable pour tout objet de $D^b(\gsp(V\otimes Q))$.
\end{remarque}

\noindent Désignons par $\mathrm{IH}^i$ le $i$-ème groupe de cohomologie d'intersection de~$\acal_V^*$.

\begin{corollaire} \label{corollaire_cohomologie_minimale} On a un isomorphisme canonique $\mathrm{Gal}(\bar{\Q}_p/\Q_p)$-équivariant
$$\mathrm{IH}^i\left( \acal_V^* \times \Spec(\bar{\Q}_p) ,\: \mathcal{F}_V(R) \right) \: = \: \mathrm{IH}^i\left( \acal_V^* \times \Spec(\bar{\F}_p) ,\: \RPsi(\mathcal{F}_V(R)) \right)$$
pour toute représentation $R$ de $\gsp(V\otimes \Q)$ et tout $i\in \Z$.
\end{corollaire}

\section{Filtration de Jordan-Hölder des cycles proches}

\subsection{Stratification de Kottwitz-Rapoport} \label{subsection_stratification_KR} Cette stratification de la fibre sp\'eciale~$\acal_V\times \Spec(\F_p)$ provient de la théorie du modèle local~\cite{Prang@GenestierNgo}. Notons $\tilde{W}_V$ le groupe de Weyl affine étendu du groupe réductif $p$-adique $\gsp(V\otimes \Q_p)$ et $W_\mathcal{D}^\mathrm{vec}$ le groupe de Weyl vectoriel du sous-groupe de Lévi du parabolique~$\mathrm{Stab}(\mathbb{V}^\bullet \otimes \F_p)$ de~$\gsp(V\otimes \F_p)$. Le groupe~$W_\mathcal{D}^\mathrm{vec}$ est naturellement un sous-groupe de~$\tilde{W}_V$ et est trivial dans le cas Iwahorique où~$\mathcal{D}=\{1,\cdots,g\}$.

Notons~$W_V$ le sous-ensemble fini des éléments $\mu_V$-admissibles de $\tilde{W}_V$~\cite{Prang@GenestierNgo}.
Ici, $\mu_V$ désigne le copoids minuscule~$(1,\ldots,1,0,\ldots,0)$ de $\gsp(V\otimes \Q_p)$, où~$1$ et $0$ sont chacun répétés $g$ fois.
Notons  $W_\mathcal{D}$ l'image de $W_V$ dans le double-quotient
$$W_\mathcal{D}^\mathrm{vec} \setminus \tilde{W}_V / \: W_\mathcal{D}^\mathrm{vec}\: .$$
En tant que sous-ensemble d'un groupe de Coxeter, $W_V$ hérite d'un ordre de Bruhat partiel~$\preceq$ et d'une fonction longueur $l:W_V\rightarrow \mathbb{N}$. Il existe une manière canonique d'en déduire un ordre de Bruhat~$\preceq$ et une fonction longueur~$l$ sur~$W_\mathcal{D}$~\cite[par.\ 8]{Minuscule@KottwitzRapoport}.
La proposition suivante, due à Kottwitz et Rapoport, résulte par exemple de~\cite{Prang@GenestierNgo}.

\begin{proposition2} Il existe une stratification canonique de~$\acal_V\times \Spec(\F_p)$ paramétrée par~$W_\mathcal{D}$. Notons $\acal_V^w$ la strate localement fermée associée à~$w\in W_\mathcal{D}$. Elle est lisse de dimension pure~$l(w)$ et son adhérence est
$$ \overline{\acal_V^w} \: = \: \coprod_{w' \preceq w} \: \acal_V^{w'}\: .$$
\end{proposition2} 

\begin{remarque2} On peut penser à la stratification de Kottwitz-Rapoport comme à la stratification par les singularités (horizontales) de~$\acal_V\times \Spec(\F_p)$. En effet, il résulte aisément de la théorie du modèle local que les cycles proches $\RPsi_V$ de~$\acal_V$ ont une restriction lisse à toute strate~$\acal_V^w$~\cite{Kottwitz@HainesNgo}.
\end{remarque2}

\noindent Soit $w\in W_\mathcal{D}$. Notons $j_V^w : \acal_V^w \hookrightarrow \acal_V\times \Spec(\F_p)$ l'immersion canonique. Le schéma~$\overline{\acal_V^w}$ est singulier et nous allons étudier cohomologiquement ses singularités grâce au foncteur de prolongement intermédiaire~$j_{V!*}^w$. Pour toute représentation~$R$ de~$\gsp(V\otimes \Q)$, le faisceau lisse décalé~$\mathcal{F}_V(R)[l(w)]$ est pervers sur le schéma lisse~$\acal_V^w$. Nous poserons dans la suite
$$\IC_V^w(R) \: = \: j_{V!*}^w \left( \mathcal{F}_V(R)[l(w)] \right)$$
qui est un faisceau pervers mixte sur~$\acal_V \times \Spec(\F_p)$. Ces faisceaux pervers sont des variantes avec coefficients des complexes d'intersection de l'adhérence des strates de Kottwitz-Rapoport.
En raisonnant comme dans la remarque~\ref{remarque_referee}, on peut définir le foncteur $R\mapsto \IC_V^w(R)$ sur $D^b(\gsp(V\otimes \Q))$.

\begin{remarque2} \label{remarque_purete} D'après le théorème de Gabber~\cite[5.4.3]{Pervers@BBD}, le lemme~\ref{poids_rep_faisceau}, la remarque~\ref{remarque_pas_pur} et la lissité des strates de Kottwitz-Rapoport, le complexe $\IC_V^w(R)$ est pur de poids~$-l(w)+a$ (resp. mixte de poids~$\leq -l(w)+a$ ou mixte de poids~$\geq -l(w)+a$) lorsque le complexe~$R$ est pur de poids~$-a$ (resp. mixte de poids~$\geq -a$ ou mixte de poids~$\leq -a$).
\end{remarque2}

\noindent Le lemme suivant est une généralisation directe de~\cite[4.1.2]{Poids@Morel} aux complexes d'intersection des strates de Kottwitz-Rapoport.

\begin{lemme2} \label{lemme_troncation_IC} Pour tout entier $a\in \Z$ et tout complexe~$R\in \mathrm{D}^b(\gsp(V\otimes \Q))$, on a $\mathrm{w}_{\leq a}\: \IC_V^w(R) = \IC_V^w(\mathrm{w}_{\geq -l(w)-a}  R)$ où $\mathrm{w}_{\geq -l(w)-a}$ est le foncteur de troncation sur \linebreak $\mathrm{D}^b(\gsp(V\otimes \Q))$ défini en~\ref{definition_troncation_rep}.
\end{lemme2}

\begin{demo}
Notons $K=\IC_V^w(R)$, $K_1=\IC_V^w(\mathrm{w}_{\geq -l(w)-a}R)$ et $K_2=\IC_V^w(\mathrm{w}_{<-l(w)-a}R)$. On~a $K=K_1\oplus K_2$ puisque $R= \mathrm{w}_{\geq -l(w)-a} R \oplus \mathrm{w}_{<-l(w)-a} R$. Il suffit de montrer que $K_1$ est dans ${}^w{D}^{\leq a}(\acal_V\times \Spec(\F_p))$ et que $K_2$ est dans ${}^w {D}^{>a}(\acal_V\times \Spec(\F_p))$. Il faut donc montrer que pour tout $n\in \Z$, les complexes $^{p}\mathcal{H}^n(K_1)$ et $^{p}\mathcal{H}^n(K_2)$ sont mixtes de poids respectifs $\leq a$ et $>a$. Or on a $^{p}\mathcal{H}^{n}(K_1)=\IC_V^w(\mathrm{H}^{n}(\mathrm{w}_{\geq -l(w)-a}R))=\IC_V^w(\mathrm{w}_{\geq -l(w)-a}\mathrm{H}^{n}(R))$ en utilisant la lissité des strates de Kottwitz-Rapoport et~\cite[3.5.2]{Poids@Morel}. Il suffit alors d'utiliser la remarque~\ref{remarque_purete}.
\end{demo}

\subsection{Décomposition des cycles proches} Il est facile de voir grâce à la théorie du modèle local que les complexes d'intersection des strates de Kottwitz-Rapoport permettent de décomposer les cycles proches dans le groupe de Grothendieck de la catégorie des faisceaux pervers de Weil. Introduisons au préalable la définition générale suivante.

\begin{definition2} \label{definition_polynome_faisceau} Soit~$X$ un schéma et~$K$ un complexe $\ell$-adique dans~$D^b_c(X)$. Pour tout polynôme $P=\sum_i a_i t^i \in \Z[t,t^{-1}]$, on pose
 $P(K)=\sum_i a_i K(i)$ où la somme est vue dans le groupe de Grothendieck $\mathrm{K}^0(D^b_c(X))$ de~$D^b_c(X)$ et~$K(i)$ désigne le $i$-ème twist de Tate.
\end{definition2}

Ce foncteur $K\mapsto P(K)$ de $\mathrm{K}^0(D^b_c(X))$ dans $\mathrm{K}^0(D^b_c(X))$ vérifie $P(Q(K))=(P\cdot Q)(K)$ pour tous $P,Q \in\Z[t,t^{-1}]$,  $P(K)(i)=P(K(i))$ pour tout~$i\in \Z$, $P(f_* K)=f_* K$ et $P(f_! K)=f_! K$ pour tout morphisme de type fini~$f:X\rightarrow Y$.



\begin{remarque2} Lorsque $X$ est de type fini sur~$\Spec(\F_p)$, on n'a pas~$\mathrm{w}_{\leq a} P(K) = P(\mathrm{w}_{\leq a} K)$ en général~\cite[3.1.1.(ii)]{Poids@Morel}.
\end{remarque2}

Notons~$\mathrm{Perv}_W(\acal_V \times \Spec(\bar{\F}_p))$ la catégorie des faisceaux pervers de Weil sur~$\acal_V \times \Spec(\bar{\F}_p)$, \textit{ie.} la catégorie des faisceaux pervers munis d'un isomorphisme de Frobenius au-dessus du morphisme de Frobenius de~$\acal_V \times \Spec(\bar{\F}_p)$.
Choisissons désormais un relèvement dans~$\mathrm{Gal}(\bar{\Q}_p/\Q_p)$ de l'automorphisme de Frobenius de~$\mathrm{Gal}(\bar{\F}_p/\F_p)$. Ce relèvement induit un foncteur d'oubli de la catégorie~$\mathrm{Perv}((\acal_V \times \Spec(\bar{\F}_p)) \times \bar{\eta})$ des faisceaux pervers munis d'une action compatible de~$\mathrm{Gal}(\bar{\Q}_p/\Q_p)$ dans la catégorie~$\mathrm{Perv}_W(\acal_V \times \Spec(\bar{\F}_p))$. En particulier, on peut voir les cycles proches comme des faisceaux pervers de Weil.
La proposition suivante, due à Haines et Görtz~\cite[par.\ 8.1]{Jordan@GoertzHaines}, se prouve facilement en utilisant la théorie du modèle local. Rappelons qu'on a noté $d_V=g(g+1)/2$ la dimension de~$\acal_V\times \Spec(\F_p)$.


\begin{proposition2} \label{proposition_Haines_Goertz} Pour tout $w\in W_\mathcal{D}$, il existe un unique polynôme $m_w\in \Z[t,t^{-1}]$ indépendant du choix du relèvement du Frobenius dans~$\mathrm{Gal}(\bar{\Q}_p/\Q_p)$ tel qu'on ait l'égalité
$$ \RPsi_{\acal_V} (\mathcal{F}_V(R)[d_V]) \: = \: \sum_{w\in W_\mathcal{D}} \: m_w \: (\IC_V^w(R))$$
dans le groupe de Grothendieck de la catégorie~$\mathrm{Perv}_W(\acal_V \times \Spec(\bar{\F}_p))$ des faisceaux pervers de Weil sur~$\acal_V \times \Spec(\bar{\F}_p)$ pour tout $R\in D^b(\gsp(V\otimes \Q))$.
\end{proposition2}

\section{Strates de Kottwitz-Rapoport au bord} La fibre spéciale~$\acal_V^*\times \Spec(\F_p)$ de la compactification minimale est munie de deux stratifications successives. Elle est d'abord munie de la stratification par les strates de bord, qui sont les variétés de Siegel~$\acal_{V'}\times \Spec(\F_p)$ indexées par~$V'\in \mathfrak{C}_V$. Puis chaque strate $$\acal_{V'}\times \Spec(\F_p)$$ est munie d'une stratification de Kottwitz-Rapoport de strates~$\acal_{V'}^{w'}$ pour $w'\in W_{\mathcal{D'}}$, où~$\mathcal{D'}$ est le type parahorique de~$\acal_{V'}$ décrit dans~\ref{subsection_compactification_minimale}. Le but de cette partie est de décrire l'interaction entre ces deux stratifications successives.

\subsection{Interactions combinatoires} \label{subsection_interaction_combinatoire} Soit $V'\in \mathfrak{C}_V$. Construisons une injection de~$W_\mathcal{D'}$ dans~$W_\mathcal{D}$. Le sous-espace~$V'$ de~$V$ détermine un sous-groupe parabolique~$\mathrm{Stab}(V'\otimes \Q)$ de~$\gsp(V\otimes \Q)$. Le sous-groupe de Lévi de ce parabolique est canoniquement isomorphe à~$\mathrm{GL}(V'\otimes \Q) \times \gsp((V'^\perp/V')\otimes \Q)$.
\`A~$V'$ est donc associé une injection
$$\gsp((V'^\perp/V')\otimes \Q) \: \hookrightarrow \: \gsp(V\otimes \Q)\: .$$
Soit~$g\in \gsp(V\otimes\Q)$ tel que~$g\cdot \mathrm{Stab}(V'\otimes \Q)\cdot g^{-1}$ contienne le tore standard de $\gsp(V\otimes \Q)$. La conjugaison par~$g$ induit une injection entre du groupe de Weyl affine étendu
de~$\gsp((V'^\perp/V')\otimes \Q)$ dans celui de~$\gsp(V\otimes \Q)$, 
ces deux groupes de Weyl étant défini par le choix des tores diagonaux standards.
Cette injection entre groupes de Weyl affines envoie éléments de translation sur éléments de translation.
Elle envoie même l'élément de translation~$\mu_{V'}$ (défini comme dans le paragraphe~\ref{subsection_stratification_KR} en remplaçant~$V$ par le module symplectique~$V'^\perp/V'$) sur l'élément de translation associé à un conjugué de~$\mu_V$ par le groupe de Weyl vectoriel de~$\gsp(V\otimes\Q)$.
D'après le théorème principal de~\cite{Minuscule@KottwitzRapoport}, elle induit donc une injection
$$\phi_{V'} \: : \: W_{\mathcal{D'}} \: \hookrightarrow \: W_\mathcal{D}$$
entre les sous-ensembles formés des éléments admissibles de ces groupes de Weyl.

\subsection{Interactions géométriques} Commençons par décrire une stratification de Kottwitz-Rapoport sur~$\mcalb_{V'}$ pour tout~$V'\in \mathfrak{C}_V$.

\begin{definition2} \label{definition_KR_mcalb} La strate~$\mcalb_{V'}^{w'}$ de~$\mcalb_{V'}\times \Spec(\F_p)$ paramétrée par~$w'\in W_\mathcal{D'}$ est l'antécédant de la strate~$\acal_{V'}^{w'}$ de~$\acal_{V'}\times \Spec(\F_p)$ par le morphisme structural de~$\mcalb_{V'}$ dans~$\acal_{V'}$.
\end{definition2}

La restriction de cette stratification à~$\mcal_{V'}$ admet une interprétation modulaire. En effet, puisque~$\mcal_{V'}$ est un espace de modules de $1$-motifs~\cite[par.\ 1.4]{Compact@Stroh}, il y existe un drapeau universel $H_1\:  \subset \: H_2 \: \subset \: \cdots \: \subset \: H_s\: \subset \: M[p]$ de sous-groupes finis et plats du groupe des points de $p$-torsion~$M[p]$ du $1$-motif universel~$M$.
Rappelons~\cite[par.\ 1.2.4]{Compact@Stroh} que $M[p]$ est canoniquement muni d'une filtration à trois crans, dont les gradués sont appelés partie étale, abélienne et multiplicative. Chacun des sous-groupes~$H_i$ provient d'un sous-groupe d'un de ces gradués. Il y a donc un sens à parler de la partie étale, abélienne ou multiplicative du drapeau~$H_\bullet$.
D'après~\cite[par.\ 1.2.6]{Compact@Stroh}, la position de~$V'$ dans~$V$ détermine la partie étale et la partie multiplicative de~$H_\bullet$. Par définition, 
la partie abélienne de~$H_\bullet$ provient de la chaîne universelle de sous-groupes sur~$\acal_{V'}$.

On obtient finalement une description modulaire de~$\mcalb_{V'}^{w'}\cap \mcal_{V'}$ comme espace de modules de drapeaux dans~$M[p]$ dont les parties étales et multiplicatives sont régies par~$V'$ et dont la partie abélienne est donnée par l'interprétation modulaire de la strate de Kottwitz-Rapoport~$\acal_{V'}^{w'}$ de~$\acal_{V'}\times \Spec(\F_p)$. Un facile \og bookkeeping \fg de~\cite{Compact@Stroh} permet de démontrer la proposition suivante. Avant de l'énoncer, rappelons que les hensélisés~$\acalb_V^h$ et~$\mcalb_{V'}^h$ de~$\acalb_V$ et de~$\mcalb_{V'}$ le long de certaines strates sont canoniquement isomorphes d'après la partie~\ref{subsection_compact_toro}. Les images~$\acal_V^h$ et~$\mcal_{V'}^h$ de~$\acal_V$ et~$\mcal_{V'}$ dans~$\acalb_V^h$ et~$\mcalb_{V'}^h$ se correspondent par cet isomorphisme entre~$\acalb_V^h$ et~$\mcalb_{V'}^h$.

\begin{proposition2} \label{proposition_interaction_toroidal} Soient $V'\in \mathfrak{C}_V$, $w'\in W_\mathcal{D'}$ et $w\in W_\mathcal{D}$. 
L'image dans~$\mcal_{V'}^h\times \Spec(\F_p)$ de la strate de Kottwitz-Rapoport $${\acal_{V}^w}$$ de~$\acal_V\times \Spec(\F_p)$ est vide si~$w$ n'est pas dans l'image de~$\phi_{V'}$ et
canoniquement isomorphe à la restriction à~$\mcal_V^h$ de~$\mcalb_{V'}^{w'}$ si~$w=\phi_{V'}(w')$.
\end{proposition2}

Ainsi, l'injection~$\phi_{V'}$ permet de décrire l'interaction entre strates de bord et strates de Kottwitz-Rapoport au niveau des compactifications toroïdales. De plus, dans un voisinage étale de la strate de~$\acalb_V$ paramétrée par~$V'$, la strate $$\acal_V^{\phi_{V'}(w')}$$ est un produit du tore~$E_{V'}$ décrit dans la partie~\ref{subsection_compact_toro}, d'une variété abélienne décrite dans la partie~\ref{subsection_compact_toro} et de la strate~$\acal_{V'}^{w'}$. On en déduit en considérant les dimensions une preuve géométrique du lemme suivant, qu'on pourrait de toute manière démontrer facilement en restant dans le cadre de la théorie des groupes.

\begin{lemme2} \label{lemme_longueur_dimension} Soient $V'\in \mathfrak{C}_V$ et $w'\in W_\mathcal{D'}$. Notons $w=\phi_{V'}(w')$ qui est dans $W_\mathcal{D}$. On~a
$l(w)-l(w') = d_V-d_{V'}$ où les fonctions longueurs sont relatives respectivement à~$W_\mathcal{D}$ et à~$W_\mathcal{D'}$ et où~$d_V$ et~$d_{V'}$ sont les dimension relatives de~$\acal_V$ et de~$\acal_{V'}$ sur~$\Spec(\Z[1/n])$.
\end{lemme2}

\noindent Montrons que~$\phi_{V'}$ permet également de décrire l'interaction entre strates de bord et de Kottwitz-Rapoport au niveau de la compactification minimale.

\begin{proposition2} \label{proposition_interaction_stratification} Soient $V'\in \mathfrak{C}_V$, $w'\in W_\mathcal{D'}$ et $w\in W_\mathcal{D}$. 
L'intersection dans~$\acal_V^*\times \Spec(\F_p)$ de la strate de bord~$\acal_{V'}\times \Spec(\F_p)$ et de l'adhérence~$\overline{\acal_{V}^w}$ de la strate de Kottwitz-Rapoport~$\acal_V^w$ est vide si~$w$ n'est pas dans l'image de~$\phi_{V'}$ et égale à l'adhérence $$\overline{\acal_{V'}^{w'}}$$ de~$\acal_{V'}^{w'}$ dans $\acal_{V'}\times \Spec(\F_p)$ si $w=\phi_{V'}(w')$.
\end{proposition2}

\begin{demo} Cela résulte immédiatement de la proposition~\ref{proposition_interaction_toroidal} et de la compatibilité du morphisme~$\pi:\acalb_V\rightarrow \acal_V^*$ aux stratifications du bord.
\end{demo}

\begin{remarque2} \label{remarque_strate_KR_minimale} On peut reformuler la proposition~\ref{proposition_interaction_stratification} en disant que la strate de Kottwitz-Rapoport paramétrée par~$w$ ne \og rencontre \fg pas la strate de bord~$\acal_{V'}$ si~$w$ n'est pas dans l'image de~$\phi_{V'}$ et la \og rencontre \fg exactement en la strate de Kottwitz-Rapoport paramétrée par~$w'$ si~$w=\phi_{V'}(w')$. La proposition~\ref{proposition_interaction_stratification} permet donc de définit une stratification de Kottwitz-Rapoport sur la compactification minimale~$\acal_V^*\times \Spec(\F_p)$ (qui est loin de se plonger dans une variété de drapeaux affine !) en posant
$$\acal_V^{*\: w} \: = \: \coprod_{V' \in \mathfrak{C}_V} \acal_{V'}^{\phi_{V'}^{-1}(w)}$$
où par convention $\acal_{V'}^{\phi_{V'}^{-1}(w)}$ est vide si $w\notin \mathrm{Im}(\phi_{V'})$. Le sous-schéma~$\acal_V^{*\: w}$ de~$\acal_V^*$ est une compactification canonique de la strate de Kottwitz-Rapoport~$\acal_V^w$ et le bord est une union de strates de Kottwitz-Rapoport relatives à des variétés de Siegel de genre plus petit.
\end{remarque2}

\section{Conservation des multiplicités au bord} \label{section_conservation_multiplicite}

Soient $V'\in \mathfrak{C}_V$ et $w'\in W_\mathcal{D'}$. Posons $w=\phi_{V'}(w')\in W_\mathcal{D}$.
Dans la proposition~\ref{proposition_Haines_Goertz}, nous avons défini un polynôme $m_w\in \Z[t,t^{-1}]$ relativement à~$\acal_V$ et à~$w\in W_\mathcal{D}$. On peut définir de même un polynôme $m_{w'}\in \Z[t,t^{-1}]$ relatif à~$\acal_{V'}$ et à~$w'\in W_\mathcal{D'}$.
La proposition suivante montre que ces deux polynômes sont égaux.

\begin{proposition} \label{proposition_conservation_multiplicite} Pour tout $w'\in W_\mathcal{D'}$ d'image $w\in W_\mathcal{D}$, on a $m_w=m_{w'}$ dans~$\Z[t,t^{-1}]$.
\end{proposition}

\begin{demo}
Nous allons donner une démonstration géométrique utilisant les compactifications toroïdales.
D'après la proposition~\ref{proposition_interaction_toroidal}, la strate de Kottwitz-Rapoport $\acal_V^w$ de~$\acal_V\times \Spec(\F_p)$ rencontre un voisinage étale~$U$ de la strate de bord de~$\acalb_{V}\times \Spec(\F_p)$ paramétrée par~$V'$. D'après les résultats de la partie~\ref{subsection_compact_toro}, $U$ est isomorphe à un ouvert étale de~$\mcalb_{V'}$. Les restrictions à~$U$ des cycles proches~$\RPsi_{\acal_V}$ et~$\RPsi_{\mcal_{V'}}$ se correspondent via cet isomorphisme.
D'après le théorème de changement de base lisse, on a $\RPsi_{\mcal_{V'}}(\Q_\ell)=(\pi')^*\RPsi_{\acal_{V'}}(\Q_\ell)$ où~$\pi'$ désigne le morphisme lisse 
de~$\mcal_{V'}$ dans~$\acal_{V'}$

D'après la proposition~\ref{proposition_interaction_toroidal}, les strates de Kottwitz-Rapoport paramétrées par~$w$ et~$w'$ se correspondent 
par l'isomorphisme local entre~$\acal_V$ et~$\mcal_{V'}$. De plus, la strate de Kottwitz-Rapoport de~$\mcal_{V'}$ paramétrée par~$w'$ est par définition l'image inverse par~$\pi'$ de la strate de Kottwitz-Rapoport de~$\acal_{V'}$ paramétrée par~$w'$. Ainsi, les complexes~$\IC_V^w(\mathrm{Triv})$ et 
$$(\pi')^* \: \IC_{V'}^{w'}(\mathrm{Triv})[l(w)-l(w')]$$
se correspondent par l'isomorphisme local précédent entre~$\acal_V$ et~$\mcal_{V'}$ si~$\mathrm{Triv}$ désigne la représentation triviale de~$\gsp(V\otimes \Q)$ et de~$\gsp((V'^{\perp}/V')\otimes \Q)$.
On conclut en appliquant le lemme~\ref{lemme_longueur_dimension}, puisque~$m_w$ code la multiplicité de~$\IC_V^w(\mathrm{Triv})$ et de ses twists de Tate dans la suite de Jordan-Hölder de~$\RPsi_{\acal_V}[d_V]$ et que~$m_{w'}$ code la multiplicité de~$\IC_{V'}^{w'}(\mathrm{Triv})$ dans la suite de Jordan-Hölder de~$\RPsi_{\acal_{V'}}[d_{V'}]$.
\end{demo}

\begin{remarque} Il y a donc un lien entre décomposer $\RPsi_{\acal_V}$ selon les~$\IC_V^w$ avec~$w\in W_\mathcal{D}$ et décomposer $\RPsi_{\acal_{V'}}$ selon les~$\IC_{V'}^{w'}$ avec~$w'\in W_\mathcal{D'}$.
\end{remarque}

\section{Le théorème de Pink pour les complexes d'intersection}

Rappelons qu'on a noté $j:\acal_V\hookrightarrow \acal_V^*$ l'immersion ouverte canonique. Fixons~$V'\in \mathfrak{C}_V$. On a noté $i_{V'}:\acal_{V'}\hookrightarrow \acal_V^*$ l'immersion localement fermée de la strate de~$\acal_V^*$ paramétrée par~$V'$. Dans cette partie, nous travaillons systématiquement en restriction à~$\Spec(\F_p)$, même si nous l'omettons des notations. Si~$G$ est un groupe discret et~$M$ un~$G$-module, nous noterons $\mathrm{RInv}(G,M)$ le complexe habituel de cohomologie des groupes. De même lorsque~$G$ est une algèbre de Lie.
On associe à~$V'$ le sous-groupe parabolique~$\mathrm{Stab}(V'\otimes \Q)$ de~$\gsp(V\otimes \Q)$. Notons~$N_{V'}$ le radical unipotent de ce sous-groupe parabolique. Rappelons que le groupe discret
$$\Gamma_{V'}^l$$
a été défini dans la partie~\ref{subsection_compact_toro}.
Si $R\in D^b(\gsp(V\otimes \Q))$, l'action de $\gsp(V\otimes \Q)$ sur~$R$ induit une action de $\gsp((V'^{\perp}/V')\otimes \Q)$ sur le complexe 
$$\mathrm{RInv}\left( \Gamma_{V'}^l , \:\mathrm{RInv}( \mathrm{Lie}(N_{V'}), \: R)\right)\: $$
Nous pouvons parler du complexe $\ell$-adique sur~$\acal_{V'}$ associé à cet objet de \linebreak $D^b(\gsp((V'^{\perp}/V')\otimes \Q))$ par le foncteur~$\mathcal{F}_{V'}$.

\begin{theoreme} \label{theoreme_Pink} Pour tout $R\in D^b(\gsp(V\otimes \Q))$ et tout $w'\in W_\mathcal{D'}$ d'image $w\in W_\mathcal{D}$, il existe un isomorphisme canonique
$$i_{V'}^* \circ j_* \: \IC_V^w (R)[d_{V'}-d_V] \: \isolong \:  \IC_{V'}^{w'}\left( \mathrm{RInv}\left(  \Gamma_{V'}^l , \:\mathrm{RInv}(\mathrm{Lie}(N_{V'}), \: R)\right) \right)\: .$$
Pour tout $w\in \mathcal{D}$ qui n'est pas dans l'image de~$\phi_{V'}$, on a $i_{V'}^* \circ j_* \: \IC_V^w (R)=0$.
\end{theoreme}

\begin{demo}
L'annulation lorsque~$w\notin \mathrm{Im}(\phi_{V'})$ résulte immédiatement de~\ref{proposition_interaction_stratification}.
Introduisons une compactification toroïdale~$\acalb_V$ et notons~$Z_{V'}$ sa strate paramétrée par~$V'$. On dispose d'un diagramme commutatif à carré cartésien
$$\xymatrix{Z_{V'} \ar[d]^{\pi_{V'}} \ar[r]^{I_{V'}} & \acalb_V \ar[d]^\pi  & \\
\acal_{V'} \ar[r]_{i_{V'}} & \acal_V^* & \acal_V \ar[l]^{j} \ar[lu]_J
}$$
Le théorème de changement de base propre pour~$\pi$ induit un isomorphisme canonique
\begin{equation} \label{equation_cbp}
i_{V'}^* \circ j_* \left( \IC_V^w(R)\right)  \: \isolong \: \pi_{V' *} \circ I_{V'}^* \circ J_* \left( \IC_V^w(R)\right) \: .
\end{equation}
Montrons que 
si $\mathrm{Triv}$ désigne la représentation triviale de~$\gsp(V\otimes \Q)$ et de \linebreak $\gsp((V'^{\perp}/V')\otimes \Q))$, on a un isomorphisme canonique
\begin{equation} \label{equation_bord_Pink}
I_{V'}^* \circ J_* \left( \IC_V^w(R)\right)\: \isolong \: \left( \pi_{V'}^* \: \IC_{V'}^{w'}(\mathrm{Triv})[l(w)-l(w')] \right) \otimes \left( I_{V'}^* \circ J_*\:  \mathcal{F}_V(R) \right) \: .
\end{equation}
En effet, on peut remplacer $\acalb_V$ par son hensélisé en~$Z_{V'}$. On peut donc remplacer l'immersion $\acal_V \hookrightarrow \acalb_{V}$ par l'immersion quotient
$$ \mcal_{V'}  / \: \Gamma_{V'}^{l} \: \hookrightarrow \:  \mcalb_{V'}  / \: \Gamma_{V'}^{l}\: .$$
Pour simplifier, on notera encore $I_{V'}$ et $J$ les morphismes correspondants pour~$\mcalb_{V'}/\: \Gamma_{V'}^{l}$.
Soit~$\pi'$ le morphisme lisse canonique de $\left( \mcalb_{V'} \right) / \: \Gamma_{V'}^{l}$ vers~$\acal_{V'}$. Sur $\left( \mcal_{V'} \right) /  \Gamma_{V'}^{l}$ on a
$$\IC_V^w(\mathrm{Triv})\: \isolong \: (\pi')^{*}\: \IC_{V'}^{w'}(\mathrm{Triv})[l(w)-l(w')]$$
d'après la proposition~\ref{proposition_interaction_stratification}. On conclut l'existence de l'isomorphisme~\ref{equation_bord_Pink} de la suite d'égalités
\begin{eqnarray*}
I_{V'}^* \circ J_* \left( \IC_V^w(R)\right) & \isolong & I_{V'}^* \circ J_* \left( \IC_V^w( \mathrm{Triv}) \otimes \mathcal{F}_V(R)\right) \\
& \isolong & I_{V'}^* \circ J_* \left( (\pi')^{*}\: \IC_{V'}^{w'}(\mathrm{Triv})[l(w)-l(w')]  \otimes \mathcal{F}_V(R)\right) \\
& \isolong & \left( \pi_{V'}^{*}\: \IC_{V'}^{w'}(\mathrm{Triv})[l(w)-l(w')] \right) \otimes \left( I_{V'}^* \circ J_* \left(  \mathcal{F}_V(R)\right) \right)
\end{eqnarray*}
\noindent On en déduit finalement l'égalité
$$i_{V'}^* \circ j_* \: \IC_V^w (R) \: \isolong \: \IC_{V'}^{w'}(\mathrm{Triv})[l(w)-l(w')] \: \otimes \: \left( \pi_{V' *} \circ I_{V'}^* \circ J_* \: \mathcal{F}_V(R) \right)$$
en combinant les égalités~\ref{equation_bord_Pink} et~\ref{equation_cbp}. On conclut la démonstration du théorème en appliquant le théorème de Pink classique~\cite[3.5.1]{Baily@Pink}, qui affirme que
$$i_{V'}^* \circ j_*   \left(  \mathcal{F}_V(R)\right) \: \isolong \: \pi_{V' *} \circ I_{V'}^* \circ J_* \left(  \mathcal{F}_V(R)\right) \: \isolong \:  \mathcal{F}_{V'}\left( \mathrm{RInv}\left( \Gamma_{V'}^l , \:\mathrm{RInv}(\mathrm{Lie}(N_{V'}), \: R)\right)\right)$$
et le lemme~\ref{lemme_longueur_dimension}.
\end{demo}

\begin{remarque} Nous avons appliqué le théorème de Pink à la variété de Siegel~$\acal_V$ qui a mauvaise réduction sur~$\Spec(\F_p)$. En effet, bien que~\cite{Baily@Pink} soit formulé sur~$\Spec(\Q)$, le seul ingrédient utilisé dans la démonstration est l'existence d'une compactification minimale et de compactifications toroïdales satisfaisant certaines relations d'ordre combinatoire. Dans le cas de~$\acal_V$, ces compactifications existent sur~$\Spec(\F_p)$ d'après~\cite{Compact@Stroh} et~\cite{Minimale@Stroh} et satisfont les mêmes relations combinatoires qu'en caractéristique nulle.
\end{remarque}

\section{Prolongement intermédiaire au bord}

Commençons par expliquer comment les complexes $\IC_V^w$ se restreignent au bord de la compactification minimale lorsqu'on leur applique le foncteur de prolongement intermédiaire~$j_{!*}$. Remarquons au préalable que le foncteur $R\mapsto j_{!*}\: \IC_V^w(R)$ défini sur la catégorie des représentations de $\gsp(V\otimes \Q)$ est exact. Il se dérive donc en en foncteur triangulé et passe au groupe de Grothendieck. C'est ce morphisme défini entre groupes de Grothendieck que nous allons étudier.
Avant d'énoncer la proposition, introduisons quelques notations.

Pour tout $V'\in \mathfrak{C}_V$ notons $\mathscr{C}_{V'}$ l'ensemble des drapeaux 
$V^\bullet=( 0 \subsetneq V^r \subsetneq \cdots \subsetneq V^0)$ d'éléments de~$\mathfrak{C}_V$ telles que~$V^0$ soit conjugué à~$V'$ par le sous-groupe de niveau~$\Gamma_V$.
L'ensemble~$\mathscr{C}_{V'}$ est muni d'une action de~$\Gamma_V$ et le quotient est fini.
Tout élément~$V^\bullet$ de~$\mathscr{C}_{V'}$ détermine un sous-groupe parabolique~$P_{V^\bullet}=\mathrm{Stab}(V^\bullet\otimes\Q)$ de~$\gsp(V\otimes \Q)$. Notons $N_{V^\bullet}$ le radical unipotent de~$P_{V^\bullet}$. Le sous-groupe
de Lévi de~$P_{V^\bullet}$ est canoniquement isomorphe au produit
$$\mathrm{GL}(V^r)\times \mathrm{GL}(V^{r-1}/V^r)\times \cdots \times \mathrm{GL}(V^0/V^1)\times \gsp(V^{0\perp}/V^0)\: .$$
Notons~$\Gamma_{V^\bullet}^l$ l'intersection
de~$\Gamma_V$ avec $\mathrm{GL}(V^r)\times \mathrm{GL}(V^{r-1}/V^r)\times \cdots \times \mathrm{GL}(V^0/V^1)$.
Lorsque~$V^\bullet$ est bien déterminé modulo l'action de~$\Gamma_V$, les groupes $N_{V^\bullet}$ et~$\Gamma_{V^\bullet}^l$ sont bien déterminés à conjugaison près par~$\Gamma_V$. En particulier, la classe d'isomorphisme de leurs groupes de cohomologie est bien déterminée.
Notons enfin~$\sharp V^\bullet$ l'entier~$r$ associé au drapeau $V^\bullet=( 0 \subsetneq V^r \subsetneq \cdots \subsetneq V^0)$ et $\delta_k(V^\bullet)=\mathrm{dim}(V^k\otimes \Q)$ pour tout $1\leq k\leq r$.

Pour tout $1\leq \delta\leq g=\mathrm{dim}(V\otimes \Q)/2$, notons~$S_\delta$ le tore de~$\gsp(V\otimes \Q)$ formé des matrices scalaires
$$\mathrm{diag}(1,\cdots,1,\lambda,\cdots,\lambda,\lambda^2,\cdots,\lambda^2)$$
où~$\lambda\in \Gm$, les entrées $1$ et~$\lambda^2$ sont répétées chacunes~$\delta$ fois et l'entrée~$\lambda$ est répétée~$2g-\delta$ fois. Si $R\in D^b(\gsp(V\otimes \Q))$ et~$a\in \Z$, on note $\mathrm{w}_{<a}^\delta(R)$ le sous-espace propre généralisé de~$R$ où~$S_\delta$ agit par des caractères de la forme
$$\mathrm{diag}(1,\cdots,1,\lambda,\cdots,\lambda,\lambda^2,\cdots,\lambda^2) \: \mapsto \: \lambda^b$$
avec $b<a$. De même pour~$\mathrm{w}_{\geq a}^\delta(R)$. Les différents foncteurs triangulés~$\mathrm{w}_{<a}^\delta$ commutent entre-eux lorsque~$d$ varie. Si~$R\in D^b(\gsp(V\otimes \Q))$,~$V^\bullet\in \mathscr{C}_V$ et~$a\in \Z$, on note pour $r=\sharp V^\bullet$
$$R_{a,V^\bullet}=\mathrm{w}_{\geq d_V-d_{V^0}-a}^{\delta_0(V^\bullet)}\circ \mathrm{w}_{<d_V-d_{V^1} -a}^{\delta_1(V^\bullet)}\circ \cdots \circ \mathrm{w}_{< d_V-d_{V^{r}} -a}^{\delta_{r}(V^\bullet)}\: .$$

\begin{proposition} \label{proposition_IC_au_bord} Pour tout entier~$a\in \Z$, tout $R\in D^b(\gsp(V\otimes \Q))$ tel que~$\mathrm{H}^n(R)$ soit pur de poids~$a$ pour tout~$n\in \Z$, tout $V'\in \mathfrak{C}_V$ et tout $w\in W_\mathcal{D}$, le complexe
$$i_{V'}^* \circ j_{!*} \left( \IC_V^w(R) \right)[d_{V'}-d_V]$$
est nul si~$w$ n'est pas dans l'image de~$\phi_{V'}$ et égal à
$$\sum_{V^\bullet\in \mathscr{C}_{V'}/\Gamma_V} \: (-1)^{\sharp V^\bullet}\cdot \IC_{V'}^{w'}\left( \mathrm{RInv}\left(  \Gamma_{V^\bullet}^l , \:\mathrm{RInv} (  \mathrm{Lie}(N_{V^\bullet}), \: R)_{a,V^\bullet}\right)  \right) 
$$
dans le groupe de Grothendieck de~$D^b_m(\acal_{V'}\times \Spec(\F_p))$ lorsque~$w=\phi_{V'}(w')$ pour~$w'\in W_\mathcal{D'}$.
\end{proposition}

\begin{demo} 
L'annulation lorsque~$w\notin \mathrm{Im}(\phi_{V'})$ résulte immédiatement de~\ref{proposition_interaction_stratification}. Supposons désormais $w=\phi_{V'}(w')$ pour~$w'\in W_\mathcal{D'}$.

Soit $V^\bullet=( 0 \subsetneq V^r \subsetneq \cdots \subsetneq V^0)$ un élément de~$\mathscr{C}_V$ intervenant dans la somme.
Notons~$i_k$ l'immersion localement fermée de~$\acal_{V^k}$ dans~$\acal_V^*$ pour tout $0\leq k\leq r$.
D'après la proposition~\ref{proposition_Morel_alternee}, qui s'applique en poids~$b:=-l(w)+a$ en vertu de la remarque~\ref{remarque_purete} grâce à l'hypothèse de pureté faite sur~$R$, il suffit de prouver que la différence virtuelle~$\Delta_{V^\bullet}=A_{V^\bullet}-B_{V^\bullet}$ des complexes
$$A_{V^\bullet}\: = \: i_{0}^* \circ(i_{1*} \circ \mathrm{w}_{>b} \circ i_{1}^*)\circ \cdots \circ (i_{r*} \circ \mathrm{w}_{>b} \circ i_{r}^*) \circ j_* \left( \IC_V^w(R) \right)[d_{V'}-d_V]$$
et 
$$B_{V^\bullet} \: = \: i_{0}^* \circ(i_{0*} \circ \mathrm{w}_{>b} \circ i_{0}^*)\circ (i_{1*} \circ \mathrm{w}_{>b} \circ i_{1}^*)\circ \cdots \circ (i_{r*} \circ \mathrm{w}_{>b} \circ i_{r}^*) \circ j_* \left( \IC_V^w(R) \right)[d_{V'}-d_V]$$
vérifie
$$\Delta_{V^\bullet} \: = \: \IC_{V'}^{w'}\left( \mathrm{RInv}\left(  \Gamma_{V^\bullet}^l , \:\mathrm{RInv} (\mathrm{Lie}(N_{V^\bullet}), \: R)_{a,V^\bullet}\right)  \right) $$
dans le groupe de Grothendieck de~$D^b_m(\acal_{V^0}\times \Spec(\F_p))$.
On démontre ceci grâce à un procédé itératif comme dans~\cite[4.2.3]{Poids@Morel}. En effet, d'après le théorème~\ref{theoreme_Pink}, on a 
$$i_{r}^* \circ j_* \left( \IC_V^w(R) \right)[d_{V^r}-d_V] \: \isolong \:\IC_{V^r}^{w^r}\left( \mathrm{RInv}\left( \Gamma_{V^r}^l , \:\mathrm{RInv}(\mathrm{Lie}(N_{V^r}), \: R)\right) \right)
$$
où $w^r$ vérifie $w=\phi_{V^r}(w^r)$. D'après le lemme~\ref{lemme_troncation_IC} appliqué à la variété de Siegel~$\acal_{V^r}$, on a
\begin{eqnarray*}
\mathrm{w}_{>b}\circ i_{r}^* \circ j_* \left( \IC_V^w(R) \right)[d_{V^r}-d_V] &  \isolong & \IC_{V^r}^{w^r}\left( \mathrm{w}_{< -l(w^r)-b}^{\delta_r(V^\bullet)} \: \mathrm{RInv}\left(  \Gamma_{V^r}^l , \:\mathrm{RInv}(\mathrm{Lie}(N_{V^r}), \: R)\right) \right) \\
&  \isolong & \IC_{V^r}^{w^r}\left( \mathrm{RInv}\left(  \Gamma_{V^r}^l , \:  \mathrm{w}_{< d_V-d_{V^r}-a}^{\delta_r(V^\bullet)} \: \mathrm{RInv}(\mathrm{Lie}(N_{V^r}), \: R)\right) \right)
\end{eqnarray*}
par le lemme~\ref{lemme_longueur_dimension}.
Notons $R^r= \mathrm{RInv}\left( \Gamma_{V^r}^l , \:  \mathrm{w}_{<d_V-d_{V^r}-a}^{\delta_{r}(V^\bullet)} \: \mathrm{RInv}(\mathrm{Lie}(N_{V^r}), \: R)\right)$ qui est un objet de
$$D^b(\gsp((V^{r\perp}/V^r)\otimes \Q))\: .$$
Pour démontrer le théorème, il faut donc calculer
\begin{equation} \label{equation_demonstration}
\mathrm{w}_{>b} \circ i_{r-1}^* \circ  i_{r*}\left( \IC_{V^r}^{w^r}(R^r)\right)
\end{equation}
et itérer les étapes précédentes.
Pour calculer~\ref{equation_demonstration}, il suffit de travailler sur l'adhérence de~$\acal_{V^{r}}$ dans~$\acal_V^*$. 
Mais cette adhérence n'est autre que la compactification minimale~$\acal_{V^r}^*$ de~$\acal_{V^r}$ et~$\acal_{V^{r-1}}$ n'est autre que la strate de bord de~$\acal_{V^r}^*$ paramétrée par l'image de~$V^{r-1}$ dans~$V^{r\perp}/V^r$, qui définit un élément de~$\mathfrak{C}_{V^r}$.

(Les deux assertions précédentes ne sont pas énoncées dans~\cite{Minimale@Stroh} mais on les démontre facilement grâce aux résultats de cet article ; avec les notations de \emph{loc. cit.}, il faut remarquer que la proposition~\cite[3.6]{Minimale@Stroh} est valable avec la même démonstration lorsque~$\bar{x}$ est un point géométrique de l'adhérence de~$Z_{V'}$ dans la compactification minimale ; le corollaire~\cite[3.7]{Minimale@Stroh} et sa démonstration se généralisent alors pour donner un isomorphisme entre l'adhérence de~$Z_{V'}$ et la compactification minimale de~$\acal_{V',n,0}$ ; comparer tout cela avec~\cite[p. 152]{Deg@FaltingsChai} où le même type d'énoncé est prouvé dans le cas de bonne réduction).

Ainsi, le calcul de~\ref{equation_demonstration} équivaut au calcul de
$\mathrm{w}_{>b}\circ i_{r}^* \circ j_* \left( \IC_V^w(R) \right)$ aux changements de notation près.
On continue ainsi de suite et l'on trouve les deux isomorphismes suivants
\begin{eqnarray*} 
A_{V^\bullet} & \isolong & \IC_{V'}^{w'}\left( \mathrm{RInv}\left( \: \Gamma_{V^\bullet}^l , \: \mathrm{w}_{< d_V-d_{V^1} -a}^{\delta_1(V^\bullet)}\circ \cdots \circ \mathrm{w}_{< d_V-d_{V^r}-a}^{\delta_r(V^\bullet)}\: \mathrm{RInv} ( \mathrm{Lie}(N_{V^\bullet}), \: R) \right)  \right) \\
B_{V^\bullet} & \isolong & \IC_{V'}^{w'}\left( \mathrm{RInv}\left( \: \Gamma_{V^\bullet}^l , \: \mathrm{w}_{< d_V-d_{V^0}-a}^{\delta_0(V^\bullet)}\circ  \cdots \circ \mathrm{w}_{< d_V-d_{V^r}-a}^{\delta_r(V^\bullet)}\: \mathrm{RInv} ( \mathrm{Lie}(N_{V^\bullet}), \: R) \right)  \right)
\end{eqnarray*}
On en déduit le théorème puisque $\mathrm{w}_{\geq d_V-d_{V^0} -a}^{\delta_{0}(V^\bullet)}=\mathrm{Id}-\mathrm{w}_{< d_V-d_{V^0}-a}^{\delta_0(V^\bullet)}$.
On remarque que la réponse finale pour le calcul de $$i_{V'}^* \circ j_{!*} \left( \IC_V^w(R) \right)[d_{V'}-d_V]$$ ne fait pas intervenir de décalage puisque~$(d_{V^r}-d_{V})+(d_{V^{r-1}}-d_{V^{r}})+\cdots+(d_{V^0}-d_{V^1})=d_{V'}-d_V$.\end{demo}

\noindent \'Enoncons et montrons à présent le théorème principal de cet article. Là encore, le foncteur $R\mapsto  j_{!*} \circ \RPsi_{\acal_V} \left( \mathcal{F}_V(R) \right)$ défini sur les représentations de~$\gsp(V\otimes \Q)$ est exact, se dérive en un foncteur triangulé puis passe au groupe de Grothendieck.

\begin{theoreme} \label{theoreme_principal} Pour tout entier~$a\in \Z$, tout $R\in D^b(\gsp(V\otimes \Q))$ tel que~$\mathrm{H}^n(R)$ soit pur de poids~$a$ pour tout~$n\in \Z$, et tout $V'\in \mathfrak{C}_V$, le complexe
$$ i_{V'}^* \circ j_{!*} \circ \RPsi_{\acal_V} \left( \mathcal{F}_V(R) \right)$$
est égal  à
$$\sum_{V^\bullet\in \mathscr{C}_{V'}/\Gamma_V} \: (-1)^{\sharp V^\bullet}\cdot \RPsi_{\acal_{V'}}\left( \mathcal{F}_{V'}\left( \mathrm{RInv}\left( \: \Gamma_{V^\bullet}^l , \:\mathrm{RInv} ( \mathrm{Lie}(N_{V^\bullet}), \: R)_{a,V^\bullet}\right)  \right) \right)
$$
dans le groupe de Grothendieck de la catégorie des faisceaux pervers de Weil sur~$\acal_{V'}\times \Spec(\bar{\F}_p)$.
\end{theoreme}

\begin{demo} D'après la proposition~\ref{proposition_Haines_Goertz}, on a
$$\RPsi_{\acal_V} (\mathcal{F}_V(R)[d_V]) \: = \: \sum_{w\in W_\mathcal{D}} \: m_w \: (\IC_V^w(R))\: .$$
Puis que $i_{V'}^*\circ j_{!*}$ commute à~$m_w$, on en déduit
\begin{equation} \label{equation_demo_finale_1}
i_{V'}^*\circ j_{!*} \circ \RPsi_{\acal_V} (\mathcal{F}_V(R)[d_V]) \: = \: \sum_{w\in W_\mathcal{D}} \: m_w \: (i_{V'}^*\circ j_{!*}\: \IC_V^w(R))\: .
\end{equation}
On déduit de la proposition~\ref{proposition_IC_au_bord} que~\ref{equation_demo_finale_1} est égal à
$$\sum_{w'\in W_\mathcal{D'}} \: m_{\phi_{V'}(w)} \left( \sum_{V^\bullet\in \mathscr{C}_{V'}/\Gamma_V} \: (-1)^{\sharp V^\bullet+d_V-d_{V'}}\cdot \IC_{V'}^{w'}\left( \mathrm{RInv}\left( \Gamma_{V^\bullet}^l , \:\mathrm{RInv} (\mathrm{Lie}(N_{V^\bullet}), \: R)_{a,V^\bullet}\right)  \right) \right) $$
donc égal à 
$$ \sum_{V^\bullet\in \mathscr{C}_{V'}/\Gamma_V} \: (-1)^{\sharp V^\bullet+d_V-d_{V'}} \cdot \sum_{w'\in W_\mathcal{D'}} \: m_{\phi_{V'}(w)} \left( \IC_{V'}^{w'}\left( \mathrm{RInv}\left(  \Gamma_{V^\bullet}^l , \:\mathrm{RInv} (\mathrm{Lie}(N_{V^\bullet}), \: R)_{a,V^\bullet}\right)  \right) \right)$$
donc égal à
$$\sum_{V^\bullet\in \mathscr{C}_{V'}/\Gamma_V} \: (-1)^{\sharp V^\bullet+d_V}\cdot \RPsi_{\acal_{V'}}\left( \mathcal{F}_{V'}\left( \mathrm{RInv}\left( \Gamma_{V^\bullet}^l , \:\mathrm{RInv} (\mathrm{Lie}(N_{V^\bullet}), \: R)_{a,V^\bullet}\right)  \right) \right)
$$
d'après la proposition~\ref{proposition_conservation_multiplicite}.
\end{demo}

\noindent Nous pouvons en déduire sans peine la trace semi-simple du Frobenius~\cite[par.\ 3.1]{Kottwitz@HainesNgo}
$$\mathrm{Tr}^{ss}_{r,V} \: : \: \acal_V^*(\F_{p^r})  \: \longrightarrow \: \bar{\Q}_\ell$$
associée au complexe ${\RPsi_{\acal_V^*}\circ j_{!*} (\mathcal{F}_V(R))}$ et à un entier~$r\geq 1$. Soit~$V'\in \mathfrak{C}_V$. Notons
$$\mathrm{Tr}_{r,V'} \: : \: \acal_{V'}(\F_{p^r}) \: \longrightarrow \: \bar{\Q}_\ell$$
la fonction trace de Frobenius habituelle associée à la somme virtuelle
$$\sum_{V^\bullet\in \mathscr{C}_{V'}/\Gamma_V} \: (-1)^{\sharp V^\bullet}\cdot  \mathcal{F}_{V'}\left( \mathrm{RInv}\left( \Gamma_{V^\bullet}^l , \:\mathrm{RInv} (\mathrm{Lie}(N_{V^\bullet}), \: R)_{a,V^\bullet}\right)  \right) $$
de faisceaux localement constant sur~$\acal_{V'}\times \Spec(\F_p)$.
 
\begin{corollaire} \label{corollaire_trace_ss} Soit $V'\in \mathfrak{C}_V$. La restriction de la fonction $\mathrm{Tr}^{ss}_{r,V}$ à la composante de bord~$\acal_{V'}(\F_{p^r})$ de~$\acal_V^*(\F_{p^r})$ est égale au produit
$$ \sqrt{p^{r \cdot d_V}} \cdot  \mathrm{Tr}_{r,V'} \cdot z_{V'}$$
où~$z_{V'}$ est la fonction de Bernstein~\cite[par.\ 4.2]{Kottwitz@HainesNgo} associée au groupe~$\gsp((V'^\perp/V')\otimes \Q)$ et au copoids minuscule~$\mu_{V'}$ défini dans la partie~\ref{subsection_stratification_KR}.
\end{corollaire}

\begin{demo} Il suffit de combiner le corollaire~\ref{corollaire_commutation_minimale}, le théorème~\ref{theoreme_principal} et le résultat principal de~\cite{Kottwitz@HainesNgo}, c'est-à-dire la  démonstration de la conjecture de Kottwitz par Haines et Ngô.
\end{demo}

La fonction $\mathrm{Tr}_{r,V'}$ s'exprime très simplement en termes de théorie des groupes. En effet, si~$x\in \acal_{V'}(\F_{p^r})$ et 
~$\gamma_0$ désigne la classe de conjugaison stable de~$\gsp((V'^\perp/V')\otimes \Q)$ associée à~$x$ par Kottwitz~\cite{Counting@Kottwitz}, on a
\begin{equation} \label{equation_terme_trace}
\mathrm{Tr}_{r,V'} \: = \: \sum_{V^\bullet\in \mathscr{C}_{V'}/\Gamma_V} \: (-1)^{\sharp V^\bullet}\cdot \mathrm{Tr}\left( \gamma_0\: , \:  \mathrm{RInv}\left(  \Gamma_{V^\bullet}^l , \:\mathrm{RInv} (\mathrm{Lie}(N_{V^\bullet}), \: R)_{a,V^\bullet}\right)  \right) \: .
\end{equation}
Cette trace fait donc apparaître une somme sur les classes de conjugaison par~$\Gamma_V$ de sous-groupes paraboliques de~$\gsp(V\otimes \Q)$ définis sur~$\Spec(\Q)$ dont le bloc central est isomorphe à~$\gsp((V'^\perp/V')\otimes \Q)$. La fonction $\mathrm{Tr}_{r,V'}$ est constante sur la partition de $\acal_{V'}(\F_{p^r})$ selon les différentes classes d'isogénie des fibres de la variété abélienne universelle.

La fonction de Bernstein~$z_{V'}$ est quant à elle constante sur les différentes strates de Kottwitz-Rapoport~$\acal_{V'}^{w'}(\F_{p^r})$ de~$\acal_{V'}(\F_{p^r})$.

On voit au final que la fonction~$\mathrm{Tr}^{ss}_{r,V}$ sur~$\acal_{V}^*(\F_{p^r})$ ne dépend que de l'intersection des strates de bord~$\acal_{V'}(\F_{p^r})$ et de leur partition en classe d'isogénies avec les strates  de Kottwitz-Rapoport~$\acal_{V}^{*\: w}(\F_{p^r})$ définies dans la remarque~\ref{remarque_strate_KR_minimale}.

Le corollaire~\ref{corollaire_trace_ss} permet facilement d'obtenir une expression de la fonction~$L$ semi-simple~\cite[par.\ 9.3.1]{Clay@Haines} associée à la cohomologie d'intersection
$$\mathrm{IH}^\bullet \left(\acal_V^*\times \Spec(\bar{\Q}_p)\: ,\: \mathcal{F}_V(R)\right)$$
de la compactification minimale. Cette expression prend la forme d'une somme sur les classes de conjugaison par~$\Gamma_V$ de sous-groupes paraboliques de~$\gsp(V\otimes \Q)$ définis sur~$\Q$, puis d'une somme sur les triplets de Kottwitz~$(\gamma_0;\gamma,\delta)$ associés à un groupe symplectique du type~$\gsp((V'^\perp/V')\otimes \Q)$ qui vérifient la condition diagonale~$\alpha(\gamma_0;\gamma,\delta)=0$~\cite{Counting@Kottwitz}. Les termes que l'on somme sont des produits de termes de trace~\ref{equation_terme_trace}, d'intégrales orbitales en toutes les places finies de~$\Q$ différentes de~$p$ et d'intégrales orbitales tordues de fonctions de Bernstein~$z_{V'}$ en~$p$. 

Pour comparer cette expression à la formule des traces d'Arthur et en déduire une expression automorphe de la fonction~$L$ semi-simple de la cohomologie d'intersection, il faut la stabiliser dans l'esprit de~\cite{Livre@Morel}. Il est nécessaire pour cela de démontrer des lemmes fondamentaux tordus pour les fonctions de Bernstein.

\begin{remarque} On aurait pu traiter sans difficulté majeure le cas de correspondances de Hecke comme dans~\cite[par.\ 5]{Poids@Morel}. Il faudrait toutefois changer notre description élémentaire de la combinatoire des compactifications en une description adélique.
\end{remarque}

\begin{remarque} \label{remarque_generalisation_PEL} Aux lourdeurs de notation près, les résultats et les démonstrations de cet article sont valables pour toutes les variétés de Shimura PEL de type (A) et (C) associées à un groupe~$G$ sur~$\Q$ non ramifié en~$p$ et à un sous-groupe compact ouvert de~$G(\mathbb{A}_f)$ parahorique et inclus dans un sous-groupe hyperspécial en~$p$. En effet, ces hypothèses garantissent l'existence de compactifications toroïdales et minimales (voir l'introduction de~\cite{Minimale@Stroh}).
\end{remarque}

\providecommand{\bysame}{\leavevmode ---\ }
\providecommand{\og}{}
\providecommand{\fg}{}
\providecommand{\smfandname}{et}
\providecommand{\smfedsname}{\'eds.}
\providecommand{\smfedname}{\'ed.}
\providecommand{\smfmastersthesisname}{M\'emoire}
\providecommand{\smfphdthesisname}{Th\`ese}

\end{document}